\documentclass[11pt]{article}
\usepackage{psfig,graphics,natbib}
\usepackage[scanall]{psfrag}
\usepackage[colorlinks=true,citecolor=blue]{hyperref}
\bibpunct{(}{)}{;}{x}{,}{,}
\bibpunct{(}{)}{;}{x}{,}{,}
\usepackage{amsmath}
\usepackage{amssymb}
\usepackage{fancybox}
\usepackage{fancyheadings}
\def\theequation{\arabic{equation}}

\usepackage{color}

\newtheorem{theorem}{\em Theorem}
\newtheorem{lemma}[theorem]{\em Lemma}
\newtheorem{proposition}[theorem]{\em Proposition}
\newtheorem{corollary}[theorem]{\em Corollary}

\newcommand{\mQ}{\mathcal{Q}}
\newcommand{\mY}{\mathcal{Y}}
\newcommand{\mH}{\mathcal{H}}
\newcommand{\mP}{\mathcal{P}}
\newcommand{\mfP}{\mathbb{P}}
\newcommand{\mD}{\mathcal{D}}
\newcommand{\mK}{\mathcal{K}}
\renewcommand{\ln}{\mbox{log}}
\newcommand{\diag}{\mbox{diag}}

\newcommand{{\mbm}}{{\boldsymbol \mu}}
\newcommand{{\mbt}}{{\boldsymbol \theta}}
\newcommand{{\mbb}}{{\boldsymbol \beta}}
\newcommand{{\mbp}}{{\boldsymbol \pi}}
\newcommand{{\bS}}{{\boldsymbol \Sigma}}

\newcommand{{\mG}}{\mathcal{G}} 
\newcommand{{\mS}}{\mathcal{S}} 
\newcommand{{\mM}}{\mathcal{M}} 
\newcommand{{\wmbp}}{\widehat{\mbp}}
\newcommand{{\wQ}}{\widehat{\mQ}}
\newcommand{{\wS}}{\widehat{\bS}}

\newcommand{\bbe}{\begin{Beqnarray*}}
\newcommand{\ebe}{\end{Beqnarray*}}
\newcommand{\be}{\begin{equation}}
\newcommand{\ee}{\end{equation}}
\newcommand{\bea}{\begin{eqnarray}}
\newcommand{\eea}{\end{eqnarray}}
\newcommand{\bes}{\begin{eqnarray*}}
\newcommand{\ees}{\end{eqnarray*}}

\newcommand{\bed}{\begin{description}}
\newcommand{\eed}{\end{description}}

\begin{document}

\begin{center}
{\large \bf Model Building for Semiparametric Mixtures}

\vspace{0.1in}
Ramani S. Pilla, Francesco Bartolucci and Bruce G.\
Lindsay\footnote{Pilla is with Department of Statistics, Case Western
Reserve University, Cleveland, OH, USA (E-mail: pilla@case.edu).
Bartolucci is with Department of Economics, University of Urbino,
Urbino, Italy (E-mail: Francesco.Bartolucci@uniurb.it). Lindsay is
with Department of Statistics, Pennsylvania State University,
University Park, PA, USA (E-mail: bgl@psu.edu). Pilla's research was
partially supported by the National Science Foundation (NSF) grant DMS
02-39053 and Office of Naval Research grants N00014-02-1-0316 \&
N00014-04-1-0481. Bartolucci's research was partially supported by the
MIUR grant 2002. Lindsay's research was partially supported by the NSF
grants DMS 01-04443 and DMS 04-05637.}

\end{center}

\vspace{0.1in}
\centerline {\bf Abstract}

An important and yet difficult problem in fitting multivariate mixture
models is determining the mixture complexity. We develop theory and a
unified framework for finding the nonparametric maximum likelihood
estimator of a multivariate mixing distribution and consequently 
estimating the mixture complexity. Multivariate mixtures provide a
flexible approach to fitting high-dimensional data while offering data
reduction through the number, location and shape of the component
densities. The central principle of our method is to cast the mixture
maximization problem in the concave optimization framework with
finitely many linear inequality constraints and turn it into an
unconstrained problem using a {\em penalty function}. We establish the
existence of parameter estimators and prove the convergence properties
of the proposed algorithms. The role of a ``sieve parameter'' in
reducing the dimensionality of mixture models is demonstrated. We
derive analytical machinery for building a collection of
semiparametric mixture models, including the multivariate case, via
the sieve parameter.  The performance of the methods are shown with
applications to several data sets including the cdc15 cell-cycle yeast
microarray data.

\vspace{1ex}


\section{Introduction}
\label{sec-intro}

Multivariate mixture modeling is a bridge between clustering and
nonparametric multivariate density estimation. The estimated
multivariate mixture model provides both an estimate of the density
for the overall data and partitions the data into several components
or clusters. Determining the mixture complexity is challenging even in
one dimension
\citep{laird:78,jewell:82,titter:85,lesperance:92,roeder:94,
lindsay:95,mclach:01,pilla2:03,scott1:04,pilla:05}. A fundamental
problem in high-dimensional modeling is determining the number of
components and their centers. One popular model-free approach to
high-dimensional modeling is the {\em K-means} algorithm (see
\cite{hastie:01} and the references therein). Model-based techniques
such as density estimation \citep{scott:92,james:01,scott2:04} and
multivariate mixture models provide a reliable and flexible approach
to high-dimensional modeling while providing a data reduction through
the number, location and shape of the component densities.  In the
context of mixture models, the problem becomes determining the number
of mixture components and estimating the corresponding location
parameter vectors. Furthermore, mixture models provide much of the
flexibility of the nonparametric approaches, while retaining many
advantages of the parametric approaches
\citep{laird:78,roeder:92,lesperance:92,lindsay:95,char:05,scott2:04}. 

One of the main reasons for the popularity of model-free approaches
such as the K-means algorithm for high-dimensional modeling is the
lack of a unified and powerful technique for fitting multivariate
mixtures. The focus of this article is to develop analytical machinery
for {\em building a collection of semiparametric mixture models},
including the multivariate case. The theory and methods developed in
this article have applications to image analysis, high-dimensional
clustering and data mining, to name a few. The practical applications
of semiparametric mixture models are broad and include case-control
studies with errors-in-variables \citep{roeder:96}, random effects
models and empirical Bayes method \citep{lindsay:95}. A natural
outcome of applying a multivariate mixture model for high-dimensional
clustering is that (1) each cluster is statistically represented by a
parametric distribution; for instance, normal in continuous case and
Poisson in discrete case and (2) it provides the proportion of
observations in each cluster through the estimated mixture
probability. Furthermore, statistical tests can be developed easily
based on the parameter estimators of the multivariate mixture models
to answer various scientific or biological questions.  Due to high
levels of noise inherent in many of the massive data sets, including
the microarray technology, it is highly desirable to carry out the
high-dimensional data analysis within a statistical framework.

Let $m := |\mbox{supp}(\mQ)|$ be the size of the support set of a {\em
mixing measure} $\mQ$; i.e., the {\em mixture complexity}. If $\mQ$ is
finitely supported, $m$ is finite and otherwise it is infinite. The
current standard approach for finding the maximum likelihood estimator
(MLE) of $\mQ$, when $m$ is known a priori or fixed, is the well-known
EM algorithm developed in the seminal article by \cite{dempster:77}.

In the absence of the knowledge of mixture complexity, it is
instructive to start with an {\em overparameterized mixture model} and
search over the whole continuous parameter space effectively to obtain
a parsimonious mixture model. Overparameterization here refers to
fitting a model with many components relative to the actual number in
the nonparametric MLE (NPMLE) of $\mQ$; hence, there is a redundancy
of components in the mixture model. Such a scheme would be robust to
parameter starting values chosen for fitting the mixture model. To
accomplish this, one requires a powerful mixture algorithm that pushes
most of the mixture probabilities to zero, which is on the boundary of
the parameter space. The focus of this article is to develop theory
and create robust (to starting values) as well as powerful algorithms
to address this problem.

The popular and widely employed EM-based algorithms are particularly
slow for fitting such an overparameterized mixture model and it is
very difficult, if not impossible, to remove the unnecessary
components; since the algorithm can never reach such a boundary
point. Moreover, the EM algorithm is sensitive to parameter starting
values and fails to converge in certain mixture problems. Figure
\ref{fig1} in Section \ref{sec-poiss} demonstrates this aspect of the
algorithm.  Other examples where the EM algorithm converges to saddle
points or fails to converge are noted by \cite{mclach:97}.

\subsection{Statistical Framework}

Let $\mathcal{F} := \{f_{\boldsymbol{\theta}}({\bf x})\!:
\boldsymbol{\theta} \in \Omega \subset \Re^p\}$ be a family of
probability density functions with respect to a $\sigma$-finite
dominating measure $\mu$ for a $p$-dimensional random vector ${\bf x}
\in \mathcal{X} \subset \Re^n$ and a $p$-dimensional location
parameter vector $\boldsymbol{\theta} \in \Omega \subset \Re^p$, a
measurable space. Assume that the component density
$f_{\boldsymbol{\theta}}({\bf x})$ is bounded in $\mbt$ for each ${\bf
x} \in \mathcal{X}$. Let $\mathcal{G}$ be the space of all probability
measures on $\Omega$ with the $\sigma$-field generated by its Borel
subsets. For a given $\mQ \in \mG$, we assume that data vector ${\bf
x}$ arises from the marginal density
\bea
  \label{eq:mden}
  {\rm g}_{_\mQ}({\bf x}) = \int f_{\boldsymbol{\theta}}({\bf
  x}) \, d \mQ(\mbt) \quad \mbox{for} \quad {\bf x} \in \mathcal{X}
\eea
which is referred to as a {\em mixture density}. The mixture model
(\ref{eq:mden}) is also applicable to empirical Bayes estimation,
where $\mQ$ is an unknown prior distribution and the objective becomes
estimation of the posterior distribution of $\mbt$ without assuming a
functional form for the prior distribution.

The goal is to estimate the {\em mixing measure} $\mQ$ by finding the
probability measure $\wQ \in \mG$ that maximizes the nonparametric
mixture loglikelihood $\log \, \{L(\mQ)\} = \sum_{i = 1}^n \log \,
\{{\rm g}_{_\mQ}({\bf x}_i) \}$. It is well known that finding the
NPMLE of $\mQ$ is computationally intensive 
\citep{lesperance:92,roeder:92,lindsay:95,bickel:98,susko:99}.
Although $\mQ$ is an arbitrary probability measure, under mild
conditions \cite{lindsay1:83,lindsay2:83} showed that finding the MLE
involved a standard problem of convex optimization, that of maximizing
a concave function over a convex set. One consequence of this is that,
as long as $l(\mQ)$ is bounded, the MLE of $\mQ$ is concentrated on a
support of cardinality at most that of $d$---the {\em number of
distinct observed data vectors}. This is a very useful, albeit
surprising, result since a potentially difficult nonparametric
estimation problem is reduced to that of a finite dimension; hence
algorithms can be constructed to find the solution. Hence, we restrict
the attention to discrete probability measures $\mQ$ having
$p$-dimensional \textit{support vectors} $\mbt_{1}, \ldots, \mbt_{m}$
collected in a matrix $\boldsymbol{\Theta}$ with a corresponding
vector of \textit{masses} or \textit{mixing probabilities} denoted by
$\mbp = (\pi_{1}, \ldots, \pi_{m})^{T}$ such that $\mbp$ is in the
{\em unit simplex} ${\boldsymbol \Pi} := \{\mbp \in \Re^m:
\pi_j \in [0, 1], \sum_{j = 1}^m \pi_j = 1\}$.

We consider model fitting for both discrete and continuous
data. Therefore, it is instructive to define ${\bf
f}_{\boldsymbol{\theta}} := \{f_{\boldsymbol{\theta}}({\bf y}_{1}),
\ldots, f_{\boldsymbol{\theta}}({\bf y}_{d}) \}^{T}$ to be the
$d$-dimensional vector of distinct likelihood terms, where $T$ denotes
transpose, $({\bf y}_{1}, \ldots, {\bf y}_{d}) \in \mY$ are the
distinct observation vectors arising from the original data vectors
$({\bf x}_{1}, \ldots, {\bf x}_{n}) \in \mathcal{X}$. Let $n_{i}$ be
the number of times ${\bf y}_{i}$ occurs in the sample of ${\bf x}$
vectors. For continuous data, $n = d$ and for discrete data, often $n
\gg d$.

The observed data matrix of dimension $(n \times p)$ with 
row vectors ${\bf y}_{i}$ ($i = 1, \ldots, n$) is assumed to arise
from the mixture density ${\rm g}_{_\mQ}({\bf y}_{i}) = \sum_{j =
1}^{m} \pi_{j} \; f_{\boldsymbol{\theta}_j}({\bf y}_{i})$, and the
discrete mixing measure can be represented as $\mQ = \sum_j \pi_j \, 
\varrho(\mbt_j)$, where $\varrho(\mbt)$ is a discrete measure with
mass one at $\mbt \in \Omega$. If $m$ is fixed, the model
${\rm g}_{_\mQ}({\bf y}_{i})$ will be referred to as the
\textit{$m$-component mixture model}, and one can always find the
NPMLE of $\mQ$ using $m$ equal to $d$
\citep{lindsay1:83,lindsay2:83}. However, the actual number of
distinct support vectors with positive mixture probability, referred
to as {\em active supports}, can be as small as one. The mixture
loglikelihood of $\mQ$ becomes
\begin{equation}
  \label{eq:lleq}
  l(\mQ) = \ln \, L(\mQ) = \sum_{i = 1}^d
  n_{i} \; \ln \, \{{\rm g}_{_\mQ}({\bf y}_{i})\} \quad \mbox{over
  all} \quad \mQ \in \mathcal{G}.
\end{equation}
The goal is to find $\wQ \in \mG$ such that $l(\wQ) = \sup_{\mQ \, \in
\, \mG} \, l(\mQ)$.

The biggest practical problem one faces in solving the loglikelihood
equations in (\ref{eq:lleq}) is that the number of inequalities is
equal to the number of elements of the parameter space $\Omega$. There
are some important problems where this number is finite, although
possibly very large, such as in target recognition, hyperspectral
image analysis and positron emission tomography. For these problems
discretization of the parameter space is directly relevant. In other
problems, $\Omega$ may be a continuous space; hence, one needs a
machinery for approximating the parameter space to solve these
equations (see Section \ref{sec-parmix}). 

\subsection{Main Results}

In this article, we develop a unified framework for finding the NPMLE
of a multivariate mixing distribution $\mQ$ and consequently for
building a collection of semiparametric mixture models. The key
ingredients for building these models are the ``sieve parameter''
controlling the dimensionality of the mixture problem (as shown in
Figures \ref{fig-astrod} and \ref{fig-irisd}) and the ability to fit
overparameterized mixture models. This collection of models enable us
to investigate the role of many overlapping densities, thereby
creating an ideal situation for solving large-scale practical
problems.

We create a powerful technique referred to as the ``Penalized Dual
method'' and an efficient algorithm for fitting overparameterized
mixture models. Consequently we have a method for estimating the
mixture complexity. This algorithm is a step in the direction of
developing a unified framework for building a collection of
semiparametric mixture models.

The underlying principle for our method is that the mixture
loglikelihood forms a concave functional on the convex set of all
probability distributions which implies that there exists a {\em dual
optimization problem} [Section 5.3, \cite{lindsay:95}]. 
\cite{lesperance:92} and \cite{susko:99} exploited
this to create an elegant algorithm for finding the NPMLE of $\mQ$ in
univariate mixtures.  This research is in the same spirit but extends
these ideas by introducing a ``penalty term''.  A fundamental feature
of our approach is that it eliminates ad hoc procedures to estimate
the penalty parameter. The dual problem has a statistical
interpretation analogous to the least squares problem and the
formulation is strikingly similar to the one that arises in empirical
likelihood \citep{owen:01} framework (see Section \ref{sec-stat}).

\vspace{0.05in}
Our main results are summarized as follows.
\vspace{-0.1in}
\begin{enumerate}
\item In Section \ref{sec-opti}, we cast the mixture problem in the
  dual optimization framework. We first develop a machinery for
  approximating the continuous parameter space $\Omega$ and next
  create an algorithm (based on the Penalized Dual method) for finding
  the maximum of $l(\mQ)$. Consequently, we propose an algorithm for
  estimating the mixture complexity. We establish the convergence of
  this algorithm to MLE in Section \ref{sec-conv}. 
\vspace{-0.1in}
\item In Section \ref{sec-dual}, we develop theory for the Penalized
  Dual method in solving the dual optimization problem while
  presenting the statistical interpretation for our framework. By
  exploiting the inherent advantage of the penalty formulation, we
  derive a technique for converting parameter estimators from the
  Penalized Dual problem into the mixture probability parameters. We
  show that the Penalized Dual estimators converge to the mixture
  probability estimators as the penalty is increased, with the correct
  limits.
\vspace{-0.1in}
\item Section \ref{sec-pdalg} establishes the existence of parameter
  estimators and derives convergence results for the Penalized Dual
  algorithm, for fitting overparameterized mixture models. The
  Penalized Dual algorithm effectively yields an estimate for the
  mixture complexity. Our algorithm is based on a modification of the
  Newton-Raphson algorithm and therefore, it inherits its virtues
  while retaining the stability (i.e., monotonically increasing the
  likelihood) of EM-based algorithms. Empirical assessment of the
  faster rate of convergence of our stable and powerful algorithm
  compared to the EM algorithm will be demonstrated in Section
  \ref{sec-conv}.
\vspace{-0.1in}
\item It is shown that the algorithm based on the Penalized Dual
  method is robust to choice of parameter starting values and achieves
  the global maximum. The dimension of the dual optimization problem
  is fixed at $d$, the number of distinct observed data vectors,
  whereas for the mixture problem it grows with the mixture complexity
  $m$. For discrete mixture problems, such as binomial or Poisson, $d$
  can be much smaller than $n$ (see for example Section
  \ref{sec-poiss}). In these cases, the Penalized Dual method has no
  dimensionality cost. 
\vspace{-0.1in}
\item In Section \ref{sec-struc}, we derive several important
  structural properties of multivariate normal mixtures in which $\mQ$
  is modeled nonparametrically in the presence of an unknown
  variance-covariance matrix ${\boldsymbol \Sigma} \in \mS$ common to
  all $m$ components, where $\mS$ is a compact space. The power of our
  method rests in building a collection of semiparametric mixture
  models, including the multivariate case. We demonstrate the role of
  the sieve parameter in reducing the dimension of the mixture problem
  by creating novel graphical devices (Figures \ref{fig-astrod} and
  \ref{fig-irisd}) referred to as {\em Mixture Tree Plots}.
\vspace{-0.1in}
\item When the cardinality of the discrete parameter set (chosen for
  approximating the continuous parameter space $\Omega$) is large, the
  EM algorithm for such a mixture problem fails to converge to the
  MLE, for all practical purposes, while the Penalized Dual algorithm
  converges. From a model selection point of view, the EM algorithm
  does not eliminate the redundant components while the Penalized Dual
  algorithm yields a parsimonious mixture model. Empirical evidence
  of this is shown in Figure \ref{fig1}, Section \ref{sec-poiss}.
\vspace{-0.1in}
\item Section \ref{sec-appl} illustrates the power of the proposed
  methods using several applications. We compare our method with the
  EM algorithm, due to lack of a unified and/or stable algorithms for
  fitting the overparameterized mixture problems and for building
  semiparametric mixture models. For the univariate mixture case, we
  compared with Rotated EM algorithm (an accelerated version of the EM
  only applicable to the univariate mixtures) proposed by
  \cite{pilla:01}. Our empirical investigation demonstrate the faster
  rate of convergence of our algorithm, compared with the EM
  algorithm.
\end{enumerate}
Section \ref{sec-discuss} presents the conclusions and the Appendix
derives technical details.  

\subsection{Relevant Literature}

Widely employed model-free methods for high-dimensional modeling
include the K-means algorithm, hierarchical clustering and
agglomerative and divisive algorithms \citep{hastie:01}. However, none
of these techniques take advantage of the inherent statistical
structure of the data.

The existing model-based mixture algorithms include those for finding
the NPMLE of $\mQ$ \citep{lesperance:92,susko:99,con:01}. These
algorithms are either not fast enough for high-dimensional modeling or
not applicable for the following mixture problem: (1) the component 
densities are poorly separated and/or (2) many of the estimated
mixture probabilities are on the boundary of the parameter space. In
analyzing Sloan Digital Sky Survey data by fitting the multivariate
normal mixtures, \cite{con:01} noted that many existing techniques are
not computationally efficient and their mixture EM algorithm obtains
an improvement of only three orders of magnitude. Therefore,
developing a powerful method for fitting multivariate mixtures is
desirable. 

Although there have been some promising developments on accelerating
the EM algorithm (see \cite{mclach:97} and the references therein),
none of these methods address the overparameterized mixture problem
described earlier. To overcome the above difficulties,
\cite{pilla:96,pilla:01} proposed alternative augmentation schemes
based on the principles of the EM that provide a significantly
improved convergence rate of the EM algorithm for a class of finite
mixture models. At this time, it is not clear how to
extend these methods to multivariate mixture models; however, they do
provide an important class for comparison with our algorithm in
univariate mixture problems (comparisons are made in Section
\ref{sec-poiss}). Lindsay (1995, Section 6.3) discusses several
algorithmic methods based on directional derivatives such as the
vertex direction method and vertex exchange method to find the NPMLE
of $\mQ$. These methods also require searching over a discrete
parameter space and have certain computational disadvantages
\citep{lesperance:92}. 

\section{Mixture Maximum Likelihood Problems} 
\label{sec-opti}

In this section we first formulate the mixture problem as a convex
optimization problem and next create a framework for approximating the
continuous parameter space. Lastly, we develop an algorithm for
finding the NPMLE of $\mQ$. This algorithm forms the basis for
building a collection of semiparametric mixture models developed in
Section \ref{sec-struc}.

\subsection{Maximizing $l(\mQ)$ via Approximating $\mG$} 
\label{sec-parmix}

If the number of components in $\mQ$ is fixed, but the location
parameter vectors are unknown, then $l(\mQ)$ can have several local
maxima \citep{lesperance:92,lindsay:95,pilla:01}. Both the EM and the
K-means algorithms can get trapped at a local maximum while requiring
a priori knowledge of the mixture complexity $m$. To overcome this
problem, researchers often randomly perturb the parameter starting
values and recompute the local maxima
\citep{hall:03,hunter:04}. However, there is no theoretical 
justification to guarantee that the resulting solution reaches closer
to the global maximum. 

In fact, random parameter starting values can fail in the mixture
context for the following reasons. First, one requires an a priori
knowledge of the number of components $m$. Second, there is a danger
of choosing multiple starting values from one component while ignoring
to choose any from other components. In such a case, the EM algorithm
may not necessarily be able to locate the component from which no
parameter values are selected. This problem becomes severe when
components of unequal sizes are present; see Section \ref{sec-iris}
for an empirical investigation of this aspect.

To combat the difficulties with the parameter starting value problem
and an a priori knowledge of $m$, we develop a technique in which we
approximate $\mG$ by the set of all probability measures on a discrete
parameter space of $\Omega$. 

\vspace{0.1in}
\noindent{\sf Approximating $\mG$:} Approximate $\mG$ by $\mG_m$,
where $\mG_m$ is a set of discrete distributions generated by a finite
subset of $\Omega$. We set this finite subset to be ${\boldsymbol
\Theta}_{m} = (\mbt_1, \ldots, \mbt_m)$. As $m \to \infty$ and
${\boldsymbol \Theta}_{m}$ becomes dense in $\Omega$, the set $\mG_m
\to \mG$. In practice, a sufficiently large $m$ is chosen such that
$\mG_m$ approximates $\mG$ well. Therefore, the cardinality of
${\boldsymbol \Theta}_{m}$, namely $m$, determines how close the MLE
is to the global MLE over all measures on $\Omega$. In approximating
$\mathcal{G}$, it is important to select a suitable ${\boldsymbol
\Theta}_{m}$ while keeping computations manageable. This will be
addressed in Section \ref{sec-grid}.

\vspace{0.15in}
In what follows, we distinguish between the three mixture problems. 
\vspace{-0.05in}
\begin{enumerate}
\item The {\em fixed support mixture} problem is equivalent to maximizing
\begin{eqnarray}
  \label{eq:primal}
  l(\boldsymbol{\pi}) = \sum_{i = 1}^d n_i \; \ln \, \left\{
  {\rm g}_{_\mQ}({\bf y}_i) \right\}
\end{eqnarray}
over the parameter space ${\boldsymbol \Pi}$ while treating the support
set $\boldsymbol{\Theta}_m \subset \Omega$ as fixed. This is the
\textit{primal} or {\em mixture problem} for which we define a
``dual'' in Section \ref{sec-dual}. Note that
$\mbox{dim}(\boldsymbol{\pi}) = (m - 1)$ and grows with the
cardinality of ${\boldsymbol \Theta}_{m}$, which is a major obstacle
when $\mbox{dim}(\boldsymbol{\Theta}_m)$ is large. {\em However, the
dimension of our dual optimization problem is fixed at $d$, the number
of distinct observed data vectors}.   
\vspace{-0.1in}
\item We fix the number of components in the mixing distribution $\mQ$
to be $m$ but treat the $\mbt$ parameter vectors as unknown for
each component. Therefore, the {\em continuous support mixture model}
problem becomes simultaneously estimating $\mbp$ and $\mbt$ parameter
vectors by maximizing $l(\mQ)$ over ${\boldsymbol \Pi} \times
\boldsymbol{\Theta}_m$ for a fixed $m$. 
\vspace{-0.1in}
\item In the absence of knowledge of mixture complexity $m$,
maximizing the mixture loglikelihood in (\ref{eq:lleq}) yields an
NPMLE that is a discrete distribution on the parameter space with a
random number of component densities \citep{lindsay:95,pilla:01}. This
will be referred to as the {\em nonparametric mixture model}. The goal in turn
becomes finding the probability measure $\widehat{\mQ} \in
\mathcal{G}$ that maximizes (\ref{eq:lleq}). 
\end{enumerate}

For discrete mixture problems, such as binomial or Poisson, often $d
\ll n$; therefore, the dual methods are able to reduce the dimension
of the mixture problem. The effect of this dimensionality on the
performance of the algorithms will be demonstrated in Section
\ref{sec-poiss}.

\subsection{Characterization of the NPMLE of $\mQ$}

Let $\Gamma$ be a curve in $\Re^d$ consisting of all vectors of the
form $\{f_{\boldsymbol{\theta}}({\bf y}_1), \ldots,
f_{\boldsymbol{\theta}}({\bf y}_d)\}$, where $\mbt \in \Omega$. Under
compactness of $\Gamma$, we can define the convex hull of $\Gamma$ as
$\mbox{Conv}(\Gamma) = \{ {\bf g}_{_\mQ}\!: \mQ \in \mG, \; \mQ \;
\mbox{has finite support}\}$, where ${\bf g}_{_\mQ} = \{{\rm
g}_{_\mQ}({\bf y}_1), \ldots, {\rm g}_{_\mQ}({\bf y}_d)\}^T$. The
optimal vector ${\bf g}_{_{\wQ}} = \{{\rm g}_{_{\wQ}}({\bf y}_1),
\ldots, {\rm g}_{_{\wQ}}({\bf y}_d)\}^T \in \mbox{Conv}(\Gamma)$ and a
corresponding maximizing measure $\wQ$ can be characterized in terms
of the gradient function as shown next.

\vspace{0.1in}
\noindent{\sf Definition 1 (Finite identifiability):} For a given
family $\mathcal{F}$, suppose that $\mQ_1, \mQ_2 \in \mG$ have finite
support.  Suppose that $\mQ_j \in \mG$ yields the mixture density
${\rm g}_{_{\mQ_j}}({\bf y})$ for $j = 1, 2$. If ${\rm
g}_{_{\mQ_1}}({\bf y}) = {\rm g}_{_{\mQ_2}}({\bf y})$ for all ${\bf y}
\in \mY$ implies $\mQ_1 = \mQ_2$, then the corresponding collection of
mixture densities is said to have the {\em finite identifiability} property.

An important aspect of our technique is based on the following
fundamental property. For the NPMLE $\widehat{\mQ}$, the $i$th 
\emph{fitted model} ${\rm g}_{_{\widehat{\mQ}}}({\bf y}_{i})$ is guaranteed
to be unique (regardless of identifiability of the mixture density),
and that one can determine these fitted values by solving for the
\emph{residual} $\widehat{w}_i$, on a log-scale, defined as
\begin{equation}
  \label{eq:res}
  \ln \, \left(\widehat{w}_i \right) := \ln \, \left( \frac{n_{i}}{n}
  \right) - \ln \, \{{\rm g}_{_{\widehat{\mQ}}}({\bf y}_{i})\} \quad 
  \mbox{for} \quad {\bf y}_i \in \mY, \, i = 1, \ldots, d.
\end{equation}

In ordinary parametric likelihood problems the solution is
characterized by the likelihood equations. We extend these ideas to
our problem to show that the fitted values and the corresponding
mixing distribution $\mQ \in \mG$ can be further characterized in
terms of a set of gradient equations. That is, $\widehat{\mQ}$ is an
NPMLE if and only if
\begin{equation}
  \label{eq:psi}
  \Psi\left(\widehat{\mQ} \right) = \underset{\boldsymbol{\theta} \,
  \in \, {\boldsymbol \Theta}_m}{\mbox{sup}}
  \mD_{\mQ}(\boldsymbol{\theta}) \leq 0 \quad \mbox{for} \quad \mQ
  \in \mG,  
\end{equation}
where the \textit{gradient function}, the directional derivative of
the mixture loglikelihood in the direction of a component
density, is defined as
\begin{eqnarray}
  \label{eq:grad} 
  \mD_{\mQ}(\boldsymbol{\theta}) := \sum_{i = 1}^d n_i
  \left\{\frac{f_{\boldsymbol{\theta}}({\bf y}_i)}{{\rm g}_{_\mQ}({\bf
  y}_i)} - 1 \right\} \quad \mbox{for} \quad \boldsymbol{\theta} \in
  {\boldsymbol \Theta}_m. 
\end{eqnarray}

If a candidate maximizing measure $\widehat{\mQ}$ violates the
\textit{gradient inequality} in (\ref{eq:psi}) at some
$\boldsymbol{\theta} \in {\boldsymbol \Theta}_m$, then one is not at
the maximum. In particular, one can increase the loglikelihood by
placing some positive probability at $\boldsymbol{\theta} \in
{\boldsymbol \Theta}_m$.

\subsection{Finding the Maximum of $l(\mQ)$ via
$\boldsymbol{\Theta}_m$} 
\label{sec-alg1}

For a fixed $m$, mixture estimation is challenging due to the fact
that $l(\mQ)$ is not concave and hence there are several local maxima
\citep{lesperance:92,lindsay:95,mclach:01,pilla:01}. We create an
algorithm that is robust to the choice of parameter starting values
and reaches closer to the global maximum of $l(\mQ)$.

\vspace{0.1in}
We find the NPMLE of $\mQ$ adaptively as follows.

\vspace{0.05in}
\noindent{\sf Algorithm 1 [Finding the Maximum of $l(\mQ)$]}
\vspace{-0.1in}
\begin{enumerate}
\item Consider ${\boldsymbol \Theta}_{m} \subset \Omega$ to be the
support set of $\mQ$. Solve the fixed support mixture problem by
maximizing (\ref{eq:primal}) over the parameter space ${\boldsymbol
\Pi}$ on the support set ${\boldsymbol \Theta}_{m}$, while treating
$\mbt \in {\boldsymbol \Theta}_{m}$ as fixed. It is worth noting that
the larger the cardinality of ${\boldsymbol \Theta}_{m}$, the higher
the value of the loglikelihood at convergence.
\vspace{-0.1in}
\item Apply the MLE $\wmbp$ (with the corresponding fixed
support set $\boldsymbol{\Theta}_m \subset \Omega$) obtained in Step
1, as parameter starting values for the continuous support mixture
problem and maximize (\ref{eq:lleq}) over the product parameter space
${\boldsymbol \Pi} \times \boldsymbol{\Theta}_m$. 
\end{enumerate}
For Step 1, one requires a stable and powerful mixture algorithm and
is derived in the next sections. In particular, Algorithm 2 presented
in Section \ref{sec-algorithm2} can be employed in Step 1. The Step 2
may include estimation of other parameters in the model such as
${\boldsymbol \Sigma} \in \mS$, in the multivariate normal mixtures
context. 

For the continuous support mixture model, it will be shown in Section
\ref{sec-appl} that Algorithm 1 reaches closer to the global maximum,
if not to the global maximum. Our empirical evidence suggests that
Algorithm 1 is superior to EM-type algorithms that start with random
(or arbitrary) parameter values.

\section{The Dual Optimization Problem: Properties of Estimators}
\label{sec-dual}

We now present the problem that is dual to the primal problem
(\ref{eq:primal}) considered by \cite{lindsay1:83} and develop theory
for effectively solving it. The {\em dual problem} is to maximize
\begin{eqnarray}
   \label{eq:dual}
   l({\bf w}) = \sum_{i = 1}^d n_i \; \ln \, (w_i)
\end{eqnarray}
subject to the constraints ${\bf w} = (w_{1}, \ldots, w_{d})^{T} \in
\Re_+^{d}$ and  
\begin{eqnarray}
  \label{eq:const}
  \sum_{i = 1}^d w_i \, f_{\boldsymbol{\theta}}({\bf  y}_i)
  \leq 1 \quad \mbox{for} \quad \boldsymbol{\theta} \in \Omega.   
\end{eqnarray}

Let $\widehat{\bf w} \in \Re^d_+$ be the solution to the above dual
(or concave) optimization problem. The solution satisfies the
relationship (\ref{eq:res}), so that solving the dual problem for
$\widehat{\bf w}$ is equivalent to finding the log-scale
residuals. Hence, indirectly, via (\ref{eq:res}), we obtain the model
fitted values ${\bf g}_{_{\widehat{\mQ}}}$. A challenging step is that
one must solve for the parameter estimates for the model from these
fitted values. We create a method that exploits the particular choice
of our penalty term. Note that the constraints are linear in the
parameter vector ${\bf w} \in \Re^d_+$, and that the number of free
parameters equals $d$, while the number of constraints equals the
cardinality of ${\boldsymbol \Theta}_{m}$. The dual optimization is
with respect to ${\bf w}$ whose dimension equals $d$. This is
especially advantageous with large data sets containing, say,
thousands of observations (see Section
\ref{sec-poiss}). In Appendix A.1, we establish the relationship
between the primal and dual problems at the solution.

\subsection{Statistical Interpretation of the Dual Problem}
\label{sec-stat}

The formulation in (\ref{eq:dual}) and (\ref{eq:const}) is strikingly
similar to the one that arises in empirical likelihood framework
\citep{owen:01} in which the function $l({\bf w})$ is maximized over a
similar set of linear constraints. The empirical likelihood problem
also has a dual problem, although it does not appear to be
computationally useful. 

There is a natural interpretation of the dual problem that is
analogous to the linear model framework. In an application of the
least squares problem, one finds the fitted values $\widehat{\bf y}$
directly by projecting the data ${\bf y} \in \mY$ onto the model space
${\mathfrak X}$ (i.e., $\mfP_{\mathfrak X} \, {\bf y} = \widehat{\bf
y}$) or solves for the residual ${\bf e}$ by projecting ${\bf y}$ onto
the orthogonal complement of the model space (i.e., $\mfP_{{\mathfrak
X}^{\perp}} \, {\bf y} = {\bf e}$, where $\perp$ denotes the
orthogonal projection). In turn, we solve for the fitted values using
$\widehat{\bf y} = ({\bf y} - {\bf e})$. The primal and dual problems
have the same relationship as the projection and complementary
projection of linear models. The parallel with the linear model
framework holds if we let the data ${\bf y}_i$ equal $\ln \,
(n_{i}/n)$, the fitted model $\widehat{\bf y}_i$ equal $\ln \,
\{{\rm g}_{_{\widehat{\mQ}}}({\bf y}_{i})\}$ and the log-scale residual
$e_i$ equal $\ln \, \left(\widehat{w}_{i} \right)$. This approach
again falls very much into the spirit of the empirical likelihood,
where $(n_{i}/n)$ is the NPMLE of the probability of observing ${\bf
y}_{i} \in \mY$.

\subsection{The Penalized Dual Method: Theory}
\label{pd-theory}

The goal in this section is to turn the constrained dual optimization
problem defined in (\ref{eq:dual}) and (\ref{eq:const}) into an
unconstrained one using a ``penalty function''. This is referred
to as the \textit{Penalized Dual method}. Our method is in the spirit
of the log-barrier method \citep{renegar:00} for convex programming;
however it differs in two important respects as will be shown. 

The Penalized Dual method maximizes
\begin{eqnarray}
  \label{eq:hfunc}
   \mH_{\gamma }({\bf w}) = \sum_{i = 1}^d \left(\frac{n_i}{n} \right)
   \ln \, (w_i) - \mP({\bf  w}, \gamma)
\end{eqnarray}
over ${\bf w} \in \Re^d_+$, where $\gamma$ is a tuning parameter and 
$\mP({\bf w}, \gamma)$ is a {\em penalty function} that ensures that
the {\em Penalized Dual solution} does not violate the constraints;
the dual solution always stays in the interior of the constraint
set. One choice for the penalty function is
\begin{eqnarray}
  \label{eq:pen1}
  \mP({\bf w}, \gamma) = \frac{1}{\gamma}
  \sum_{j = 1}^m \left\{ p_{_{\boldsymbol{\theta}_j}}({\bf w})
  \right\}^{\, \gamma} \quad \mbox{for} \quad \boldsymbol{\theta}_j
  \in \boldsymbol{\Theta}_m \; \mbox{and} \;  \gamma \in \Re_+,
\end{eqnarray}
where the {\it penalty parameter} $\gamma$ is some large power and the
{\it constraint function} is defined as
\begin{eqnarray}
  \label{eq:pgrad}
  p_{_{\boldsymbol{\theta}_j}}({\bf  w}) := \sum_{i =
  1}^d w_i \, f_{\boldsymbol{\theta}_j}({\bf y}_i) > 0.
\end{eqnarray}
That is, the dual problem constraints have the form
$p_{_{\boldsymbol{\theta}_j}}({\bf w}) \leq 1$. We first show that by
increasing $\gamma$, $\mP({\bf w}, \gamma)$ will eventually create an
infinite penalty on any ${\bf w} \in \Re^d_+$ that violates the
constraints and advances the solution towards the dual problem
solution. 

The proofs for our technical results are derived in the Appendix.

\begin{proposition}
\label{lem-penalty}
{\rm For a given $\boldsymbol{\theta} \in \Omega$ and ${\bf w} \in
\Re^d_+$, the term in the summand of the penalty function $\mP({\bf
w}, \gamma)$ satisfies:  
\bes
  \frac{ \left\{ p_{_{\boldsymbol{\theta}}}({\bf w}) \right\}^{\,
  \gamma}}{\gamma} \longrightarrow \left\{ \begin{array}{ccl}  \infty
  &\mbox{if}& \quad p_{_{\boldsymbol{\theta}}}({\bf w}) > 1, \\
  0 &\mbox{if}& \quad 0 \leq p_{_{\boldsymbol{\theta}}}({\bf
  w}) \leq 1 \end{array} \right.
\ees
as $\gamma \rightarrow \infty$. When $p_{_{\boldsymbol{\theta}}}({\bf
w}) > 1$, the penalty function is increasing in $\gamma$ for $\gamma >
\{ {\rm log} \; p_{_{\boldsymbol{\theta}}}({\bf w}) \}^{-1}$. If
$p_{_{\boldsymbol{\theta}}}({\bf w}) < 1$, the penalty function is
decreasing in $\gamma$ for all $\gamma \in \Re_+$. } 
\end{proposition} 
\vspace{0.1in} 

Two main elegant features of our penalty function are the following:
(1) We can directly construct an estimator for $\pi$ parameters from
the penalized dual solution. (2) It is simple to calculate the
gradient function to assess the algorithmic convergence using the
relation (\ref{eq:gradin}), defined in Section \ref{sec-gradp}. 

It is common in the optimization literature to employ a ``barrier
function'' to build the penalty. For example, the {\em log-barrier
function} defined as
\bes
   \mP_{\star}({\bf w}, \gamma) := - \gamma
   \sum_{j = 1}^m \ln \, \{1 - p_{_{\boldsymbol{\theta}_j}}({\bf w})
   \} \quad \mbox{for} \quad j = 1, \ldots, m  
\ees
approaches $-\infty$ as ${\bf w} \in \Re^d_+$ approaches the boundary
of the feasible set from the interior \citep{roos:97,renegar:00}. The
effect of the penalty can be diminished by making $\gamma$ close to
$0$. Our focus here is on a soft penalty of the form (\ref{eq:pen1})
which is well behaved outside the feasible set; however, as will be
shown, it does force the solution into the interior.

The penalized problem is unconstrained; therefore, we can find the
``Penalized Dual optimal estimator'' denoted ${\widehat{\bf
w}}_{\gamma} = \left( \widehat{w}_{1, \gamma}, \ldots, \widehat{w}_{d,
\gamma} \right)^T$, given by (\ref{eq:fixp}) in Appendix A.2, by
solving 
\begin{eqnarray} 
  \label{eq:deriv1} 
  \frac{\partial}{\partial w_i} \mH_{\gamma}({\bf w}) =
  \frac{n_i}{n} \frac{1}{w_i} - \sum_{j = 1}^m \left\{
  p_{_{\boldsymbol{\theta_j}}} ({\bf w}) \right\}^{\, (\gamma - 1)} \,
  f_{\boldsymbol{\theta_j}} ({\bf y}_i) = 0 \quad \mbox{for} \quad i =
   1, \ldots, d.  
\end{eqnarray} 
For $\gamma = 1$, there exists an explicit solution to the above
equation as 
\begin{equation}
   \label{explicit}
   \widehat{w}_i \Big|_{\gamma = 1} = \frac{n_i}{n} \; \left\{
   \sum_{j = 1}^m f_{\boldsymbol{\theta_j}}({\bf y}_i)
   \right\}^{-1} \quad \mbox{for} \quad i = 1, \ldots, d. 
\end{equation}
This is an initial interior point solution for the algorithm. On the
other hand, the conventional log-barrier methods do not automatically
produce a starting value for the ``barrier parameter''. 

\subsection{Existence of Parameter Estimators}
\label{sec-pdes}

We consider the following method to solve for the $\pi$ parameters
from the dual problem solution by exploiting the penalized structure. 

\vspace{0.1in} 
\noindent{\sf Recovering the Primal Estimators:} Using (\ref{eq:res}),
the model fitted values ${\bf g}_{_{\wQ}}$ are found from the
penalized dual solution $\widehat{\bf w}$. However, such a solution
does not immediately provide an estimator for $\mbp \in {\boldsymbol
\Pi}$ and the technique for obtaining it is derived next.
\vspace{-0.1in}
\begin{enumerate}
\item Restrict attention to $\boldsymbol{\theta}_j \, (j = 1, \ldots,
  m)$ in ${\boldsymbol \Theta}_{m} \subset \Omega$, the support set of
  $\mQ \in \mG$, for which the constraints are tight to ensure
  $\sum_{i = 1}^d \widehat{w}_{i, \gamma} \,
  f_{\boldsymbol{\theta}_j}({\bf y}_i) = 1$.
\vspace{-0.1in} 
\item Solve for $\mbp$ using the linear equations $\sum_{j = 1}^m
\widehat{\pi}_j \, f_{\boldsymbol{\theta}_j}({\bf y}_i) =
(n_i/n)/\widehat{w}_{i, \gamma}$ for each $i = 1, \ldots, d$. 
\end{enumerate}

The penalized dual residuals are used to obtain a natural
estimator for the mixture or primal problem, denoted by 
$\widehat{\boldsymbol{\pi}}_{\gamma }^{\star} = \left(
\widehat{\pi}_{1, \gamma}^{\star}, \ldots, \widehat{\pi}_{m,
\gamma }^{\star }\right)^{T}$. The statistic, which is referred to as
the \textit{Penalized Dual estimator} is 
\begin{equation}
   \label{eq:pdest}
   \widehat{\pi}_{j,\gamma}^{\star} =
   \left\{ p_{_{\boldsymbol{\theta}_{j}}} \left(\widehat{\bf
   w}_{\gamma} \right) \right\}^{\gamma} \quad \mbox{for} \quad j = 1,
   \ldots, m, 
\end{equation}
where
\begin{eqnarray}
  \label{eq:pgrade}
  p_{_{\boldsymbol{\theta}_j}} \left(\widehat{\bf  w}_{\gamma} \right)
  = \sum_{i = 1}^d \widehat{w}_{i, \gamma} \,
  f_{\boldsymbol{\theta}_j}({\bf y}_i).  
\end{eqnarray}
In Appendix A.2, it is shown that the estimator $\widehat{w}_{i,
\gamma}$, derived in (\ref{eq:fixp}), can be approximated in terms of
$\{ p_{_{\boldsymbol{\theta}_{j}}}(\widehat{\bf w}_{\gamma}) \}^{\,
(\gamma - 1)}$. However, these latter quantities with the power
$(\gamma - 1)$ do not sum to one, and hence are turned into a {\em
candidate estimator} via normalization:
\begin{equation}
  \label{eq:ems}
  \widehat{\pi}_{j, \gamma}^{\dag} = \frac{ \left\{
  p_{_{\boldsymbol{\theta}_{j}}} \left(\widehat{\bf w}_{\gamma} \right)
  \right\}^{\,(\gamma - 1)}}{\sum_{k = 1}^{m}\, \left\{ 
  p_{_{\boldsymbol{\theta}_{k}}} \left(\widehat{\bf w}_{\gamma} \right)
  \right\}^{\,(\gamma - 1)}} \quad \mbox{for} \quad j = 1, \ldots, m.
\end{equation}

This candidate estimator is used to obtain
$\widehat{\boldsymbol{\pi}}^{\star}_{\gamma}$ with its elements having
the power $\gamma$ using the following theorem.
\begin{theorem}
\label{lemm1} 
{\rm (a) For a given $\gamma \in \Re_+$, the Penalized Dual estimator
\bes
  \widehat{\boldsymbol{\pi}}_{\gamma}^{\star} = \left[
  \left\{p_{_{\boldsymbol{\theta}_{1}}} \left(\widehat{\bf w}_{\gamma}
  \right) \right\}^{\gamma}, \ldots, \left\{
  p_{_{\boldsymbol{\theta}_{m}}} \left(\widehat{\bf w}_{\gamma} \right)
  \right\}^{\gamma} \right]^{T}
\ees
is one EM-step from the candidate estimator
$\widehat{\boldsymbol{\pi}}_{\gamma}^{\dag}$; consequently,
$\widehat{\boldsymbol{\pi}}^{\star}_{\gamma}$ yields a higher
likelihood value. (b) The estimators are in the unit simplex
${\boldsymbol \Pi}^{\star} = \{\widehat{\pi}_{j, \gamma}^{\star} \in
\Re^m: \widehat{\pi}_{j, \gamma}^{\star} \in [0, 1], \sum_{j = 1}^m \,
\widehat{\pi}_{j, \gamma}^{\star} = 1 \}$. (c) The Penalized Dual
solution $\widehat{\bf w}_{\gamma}$ satisfies
$p_{_{\boldsymbol{\theta}_{j}}}(\widehat{\bf w}_{\gamma}) \leq 1$ $(j
= 1, \ldots, m)$ and hence the estimator
$\widehat{\boldsymbol{\pi}}_{\gamma}^{\star}$ remains in the feasible
region defined by (\ref{eq:const}). } 
\end{theorem}

All the proofs are relegated to Appendix A.4. 

The estimator $\widehat{\mbp}^{\star}$ provides a direct way to obtain
the primal estimator $\widehat{\mbp}$ from our penalized dual
solution, avoiding the problems of selection and inversion.

\subsection{Properties of the Penalized Dual Estimators}
\label{sec-gradp}

In this section, we derive several statistical properties of the
estimators. First, we establish that
$\widehat{\boldsymbol{\pi}}_{\gamma}^{\star}$ converges to the MLE
$\widehat{\boldsymbol{\pi}}$ as the penalty parameter $\gamma$
increases (Theorem \ref{thm-result} below). Along the way, we
establish several important properties of the primal-gradient function
that are necessary for solving the primal-dual problem.

Although $\widehat{\boldsymbol{\pi}}_{\gamma}^{\star}$ represents an
EM improvement over $\widehat{\boldsymbol{\pi}}_{\gamma}^{\dag}$, the
candidate estimator, it is easier to establish optimization results
for the latter. The following theorem shows that for a sufficiently
large penalty, the Penalized Dual estimator will be close to the
primal estimator. Let $\widehat{\mQ}^{\dag}_{\gamma}$ be the mixing
distribution at the $\widehat{\boldsymbol{\pi}}_{\gamma }^{\dag}$
solution.

\begin{theorem} 
\label{thm-result} 
{\rm As $\gamma \rightarrow \infty$, ${\bf
g}_{_{\wQ^{\dag}_{\gamma}}} \rightarrow {\bf
g}_{_{\wQ}}$. Consequently, the candidate estimator
$\widehat{\boldsymbol{\pi}}_{\gamma}^{\dag}$ converges to the MLE
$\widehat{\boldsymbol{\pi}}$, whenever the latter is unique.}
\end{theorem} 

Our goal is to obtain the mixture estimation problem from the
penalized dual one using $\wmbp_{\gamma}^{\star}$. Therefore, it is
important to determine directly from the dual problem how accurate is
the estimator $\wmbp_{\gamma}^{\star}$. We derive the gradient 
function corresponding to the mixing distribution
$\widehat{\mQ}^{\dag}_{\gamma}$ to accomplish this. It is easier to
calculate this for $\widehat{\boldsymbol{\pi}}^{\dag}$ knowing that
$\widehat{\boldsymbol{\pi}}_{\gamma}^{\star}$ can only be better. From
the primal-gradient function in (\ref{eq:grad}), the gradient function
for the estimator $\widehat{\mQ}_{\gamma}^{\dag}$, becomes
\begin{equation}
   \label{eq:modgrad}
   \mD_{\widehat{\mQ}^{\dag}_{\gamma}}(\boldsymbol{\theta}) = \sum_{i 
   =  1}^d n_{i} \left\{ \frac{f_{\boldsymbol{\theta}}({\bf
   y}_{i})}{\sum_{k = 1}^{m}\, \widehat{\pi}_{k, \gamma}^{\dag} \, 
   f_{\boldsymbol{\theta}_{k}}({\bf y}_{i})} - 1 \right\} \quad
   \mbox{for} \quad \boldsymbol{\theta} \in {\boldsymbol \Theta}_m. 
\end{equation}

\begin{theorem}
\label{thm-grad}
{\rm The primal-gradient function at the candidate estimator 
$\widehat{\boldsymbol{\pi}}^{\dag}$ can be written as 
\begin{eqnarray}
  \label{eq:mgrad}
  \mD_{\widehat{\mQ}^{\dag}_{\gamma}}(\boldsymbol{\theta}_j) =
  \frac{p_{_{\boldsymbol{\theta}_{j}}} \left(\widehat{\bf
  w}_{\gamma} \right)}{\wp_{\gamma}} - 1 \quad \mbox{for} \quad
  \boldsymbol{\theta}_j \in {\boldsymbol \Theta}_m \, (j = 1, \ldots,
  m), 
\end{eqnarray}
where 
\begin{equation}
   \label{eq:sum}
   \wp_{\gamma}^{-1} =  \sum_{k = 1}^{m} \, \left\{
   p_{_{\boldsymbol{\theta}_{k}}} (\widehat{\bf w}_{\gamma})
   \right\}^{\,(\gamma - 1)}.  
\end{equation}
}
\end{theorem}

Theorem \ref{thm-grad} expresses the gradient function in terms of the
dual solution and leads to a simpler device for checking the accuracy
of the estimators.
\begin{corollary}
\label{cor:mgrad}
{\rm At the candidate estimator $\widehat{\boldsymbol{\pi}}^{\dag}$,
the primal-gradient function satisfies 
\begin{equation}
   \label{eq:gradin}
   \mD_{\widehat{\mQ}^{\dag}_{\gamma}}(\boldsymbol{\theta}_j) \leq
   \wp_{\gamma} - 1 \quad \mbox{for} \quad \boldsymbol{\theta}_j \in
  {\boldsymbol \Theta}_m \, (j = 1, \ldots, m),
\end{equation} 
where the term on the right-hand-side does not depend on $j$. 
}
\end{corollary}

We have established that one can refine the NPMLE of $\mQ$ to the
required accuracy by increasing $m$ and $\gamma$ appropriately.

\section{The Structure of the Penalized Dual Algorithm} 
\label{sec-pdalg}

In this section, we first investigate the structure of the penalized
dual problem viewed as a function of ${\bf w}$ and $\gamma$ and next
present a strategy for their joint estimation. Next, we present the
Penalized Dual algorithm to effectively search over the discretized
(but large) parameter space $\boldsymbol \Theta_{m} \subset
\Omega$. Lastly, convergence properties of the algorithms are derived.

We let ${\bf z} = \ln \, ({\bf w})$ to eliminate the constraint
${\bf w} \in \Re_+^{d}$.  
\begin{theorem}
  \label{lem-concave}
 (a) {\rm The function 
\begin{equation}
  \label{eq:kfunc}
  \mK({\bf z}, \gamma) = \sum_{i = 1}^d \left( \frac{n_i}{n} \right)
  \, z_i - \frac{1}{\gamma} \; \sum_{j = 1}^m \left\{
  p_{_{\boldsymbol{\theta_j}}} ({\bf z}) \right\}^{\, \gamma} \quad
  \mbox{for} \quad {\bf z} \in \Re \quad \mbox{and} \quad \gamma \in
  \Re_+ 
\end{equation}
is strictly concave in $({\bf z}, \gamma)$, where
$p_{_{\boldsymbol{\theta_j}}} ({\bf z}) = \sum_i \exp(z_i) \;
f_{_{\boldsymbol{\theta_j}}}({\bf y}_i)$. For any ${\bf z} \in \Re$ in
the feasible region defined by (\ref{eq:const}), the function
$\mK({\bf z}, \gamma)$ is strictly increasing as a function of
$\gamma$. (b) The function $\mK({\bf z}, \gamma)$ is bounded above and
achieves its maximum at ${\bf z} = \widehat{\bf z}$ and
$\gamma = \infty$.} 
\end{theorem}

\subsection{Automatic Selection of $\gamma$} 
\label{sec-pdalgor}

A fundamental aspect of our algorithm is that we can maximize
$\mK({\bf z}, \gamma)$ simultaneously with respect to ${\bf z} \in
\Re$ and $\gamma \in \Re_+$; different from the approach employed in
the conventional log-barrier methods in which this was not
possible. Therefore, we can select the penalty parameter $\gamma$
automatically. From Theorem \ref{lem-concave}, the global maximum over
${\bf z}$ and $\gamma$ is attained when ${\bf z} = \widehat{\bf z}$
and $\gamma \rightarrow \infty$.

\vspace{0.1in}
\noindent{\em Remark 1:} It may seem paradoxical to treat $\gamma$ as
an unknown parameter even though it has an optimum value of
$\infty$. When $\gamma$ is large, $\mK({\bf z}, \gamma)$ has very
severe curvature at the constraint boundary. This limits the range of
effectiveness of quadratic approximation methods.  Therefore, one
should start with a small value for $\gamma$ and increase it as the
algorithm progresses through the parameter space. This could possibly
be achieved in some other systematic fashion; however, our empirical
investigations suggest that systematic methods were not as efficient
as our approach. A possible explanation could be that our strategy
takes the curvature of the function $\mK({\bf z}, \gamma)$ into
account, in providing the relevant information for determining the
increments for $\gamma$.

\subsection{Searching Effectively Over the Discretized Parameter
Space} 
\label{sec-algorithm2}

An algorithm for efficiently searching over the large discretized
parameter space $\boldsymbol \Theta_{m} \subset \Omega$ (required for
Step 1 of Algorithm 1 described in Section \ref{sec-parmix}) is
derived next. In effect, the following algorithm is used for fitting
the fixed support mixture model. 

\vspace{0.1in} 
\noindent{\sf Algorithm 2 (The Penalized Dual Algorithm)} 
\vspace{-0.1in}
\begin{enumerate}
\item Consider $\gamma = 1$ and its corresponding explicit solution
$\widehat{\bf z}^{(1)}$ given in (\ref{explicit}) as the starting
solution for the algorithm. 
\vspace{-0.1in}
\item Maximize the concave function $\mK({\bf z}, \gamma)$
simultaneously with respect to ${\bf z} \in \Re$ and $\gamma \in
\Re_+$ using a modified Newton-Raphson algorithm [constrained the step
size to ensure monotonicity in $\mK({\bf z}, \gamma)$] until the
following convergence criterion is satisfied; namely, the $L_2$-norm
of the change in the value of $\mK({\bf z}, \gamma)$ is less than
$10^{-6}$.
\vspace{-0.1in} 
\item Fix $\gamma$ at $\gamma^{(k)}$ obtained in Step 2 and find
$\widehat{\bf z}_{\gamma^{(k)}} = \arg \, \max_{\bf z \in \Re} \; \mK
\left({\bf z}, \gamma^{(k)} \right)$ using the modified Newton-Raphson
algorithm. The algorithm is considered to have converged to the
maximum at step $t$, when the inequality based on the
primal-gradient function 
\begin{eqnarray}
  \label{eq:gstop}
  \Psi \left(\mQ^{(t)} \right) \leq 0.005 
\end{eqnarray}
is satisfied since it guarantees convergence to a similar accuracy in
the loglikelihood.
\end{enumerate}

As described in Section \ref{sec-opti}, the supremum of the gradient
function $\sup_{\boldsymbol{\theta} \, \in \, {\boldsymbol
\Theta}_m} \, \mD_{\mQ}(\boldsymbol{\theta})$ provides an assessment
of the progression to the maximum and hence the criterion in
(\ref{eq:gstop}) has a solid theoretical justification.

\vspace{0.1in}
\noindent{\em Remark 2:} The Step 3 of the algorithm is necessary
since after Step 2, the Primal Dual estimator
$\widehat{\boldsymbol{\pi}}_{\gamma}^{\star}$ obtained via 
(\ref{eq:pdest}) are often not sufficiently close to the primal
estimator $\widehat{\boldsymbol{\pi}}$. This is because the algorithm
does not necessarily satisfy the condition $\{\partial \mK({\bf z},
\gamma)/\partial {\bf z} \} = 0$ with sufficient accuracy. In our
applications, however, the primal-gradient inequality (\ref{eq:gstop})
was always achieved at the tolerance of $0.005$; in fact, often
reached significantly greater accuracy in the Penalized Dual
estimators.

\vspace{0.1in} 
\noindent{\sf The Penalized Dual Algorithm with Inactive Constraints:}
In Algorithm 2, if an estimated mixture probability $\widehat{\pi}_j
\, (j = 1, \ldots, m)$ is zero, then the corresponding constraint in
the dual problem is inactive. We can dynamically update the active
constraints by removing the inactive ones while adding new ones,
whenever the support set violated the gradient inequality. From the
Penalized Dual estimator in (\ref{eq:pdest}), it follows that if
$\widehat{\pi}_j \rightarrow 0$, then
$\{p_{_{\boldsymbol{\theta}_{j}}}(\widehat{\bf w}_{\gamma})\}^{\,
\gamma} \rightarrow 0$ which occurs when
$p_{_{\boldsymbol{\theta}_{j}}}(\widehat{\bf w}_{\gamma}) \rightarrow
0$ or $\gamma \rightarrow \infty$. In the former, one can essentially
remove the corresponding density $f_{\mbt_j}({\bf y}_i)$. It will be
shown in Section \ref{sec-appl} that the above algorithm, denoted by
PD$^{\rm IC}$, produced a further reduction in computational time. 

As a consequence of the concavity of $\mK({\bf z}, \gamma)$
established in Theorem \ref{lem-concave}, the Hessian ${\bf H}$ for
$\mK({\bf z}, \gamma)$ (derived in equation (\ref{eq:hessian}) in
Appendix A.3) is always non-singular and the sequence obtained from
the Penalized Dual algorithm (i.e., Algorithm 2) are well
defined. Even for large-scale problems such as the yeast microarray
data considered in Section \ref{sec-iris} in which ${\bf H}$ is of
dimension $697$, our modified Newton-Raphson algorithm was stable and
efficient. 

Theoretically, owing to Theorem \ref{lem-concave}, the Step 2 of the
Penalized Dual algorithm produces a sequence $\left\{{\bf z}^{(k)},
\gamma^{(k)} \right\}_{k \geq 1}$ such that the sequence of functions
$\left\{\mK({\bf z}^{(k)}, \gamma^{(k)}) \right\}_{k \geq 1} \to
\mK(\widehat{\bf z}, \infty)$ as $k \to \infty$. This effectively
implies that the sequence $\left\{{\bf z}^{(k)} \right\}_{k \geq 1}
\to \widehat{\bf z}$ and the sequence $\left\{\gamma^{(k)} \right\}_{k
\geq 1} \to \infty$ as $k \to \infty$. However, in practice,
convergence of $\gamma$ is slow; therefore, we terminate the modified
Newton-Raphson algorithm in Step 2 when $\gamma^{(k)}$ is sufficiently
large and maximize $\mK \left({\bf z}, \gamma^{(k)} \right)$ over
${\bf z}$ for a fixed $\gamma^{(k)}$.

In our experience, a direct maximization of the mixture loglikelihood
$l(\boldsymbol{\pi})$ over ${\boldmath \Pi}$ using a modified
Newton-Raphson algorithm was unstable and failed to converge to the
maximum $\wmbp$.  

\section{Convergence Properties of the Algorithms}
\label{sec-conv}

In this section, we establish the convergence properties, including
the rate of convergence, of the algorithms. 

First, we consider the algorithm for fitting the continuous support
mixture model; i.e., an algorithm employed in Step 2 of Algorithm
1. We prove that the sequence of estimates $\{ \mbb^{(s)} \}_{s \geq
1}$ obtained from the Step 2 of Algorithm 1 converges to an MLE of
$\mbb \in \Omega$, namely $\widehat{\mbb}$, for a given data ${\bf y}
\in \mY$ as $s$ increases. For instance, in the multivariate normal
mixture framework, $\mbb$ becomes $(\mQ, \boldsymbol{\Sigma})$. Assume
that the sequence of estimates $\{ \mbb^{(s)} \}_{s \geq 1}$
monotonically increases the loglikelihood $l(\mbb)$.  An algorithm is
said to converge if $\mbb^{\star} = \lim_s \mbb^{(s)}$ exists, for a
parameter vector $\mbb \in \Omega$. 

\cite{wu:83} established that monotonicity of $l(\cdot)$ does not
imply the convergence of the sequence to a stationary point; however,
if the sequence $\{l(\mbb^{(s)})\}_{s \geq 1}$ is bounded above, then
it does converge monotonically to a stationary point of $l(\mbb)$. The
convergence of $\mbb^{(s)}$ to $\widehat{\mbb}$ implies the
convergence of $l(\mbb^{(s)})$ to $l(\widehat{\mbb})$ according to the
Theorem 5, under the regularity conditions, derived by \cite{wu:83}.

Owing to Theorem \ref{lem-concave}, the Algorithm 2 (or
Step 1 of Algorithm 1) produces a sequence of estimates
$\{\mbp^{(t)}\}_{t \geq 1}$ that is guaranteed to converge to the
unique MLE $\wmbp$. Combined this result with Theorem 5 in
\cite{wu:83} establishes the convergence of Algorithm 1 to 
$\widehat{\mbb}$.

\subsection{Convergence Criteria} 
\label{sec-criteria}

For the applications and simulation experiment, we used the
convergence criterion based on the gradient function for the Penalized
Dual (PD) and discrete EM (i.e., for fitting the fixed support mixture
model) algorithms. That is, the algorithm has converged to the MLE
$\wmbp$ if the criterion in (\ref{eq:gstop}) is satisfied. For the
rest of the article, we denote the discrete EM by D-EM algorithm.

The D-EM algorithm is a sublinearly convergent algorithm
\citep{pilla:01}; therefore, a conventional convergence criterion
based on the loglikelihood change or changes in parameters, such as
\begin{eqnarray}  
  \label{eq:lstop}
  \xi^{(t)} = \big| l \left(\boldsymbol{\pi}^{(t)} \right) - l
  \left(\boldsymbol{\pi}^{(t - 1)} \right) \big| \leq \tau
\end{eqnarray}
for a given tolerance $\tau$ can be very misleading in the sense that
the actual distance to the final loglikelihood 
\bea
  \label{delta}
  \Lambda^{(t)} = \big|l\left(\widehat{\boldsymbol{\pi}} \right) -
  l\left(\boldsymbol{\pi}^{(t)} \right)\big|
\eea
can be orders of magnitude different from $\tau$. That is, this
criterion may be met even though the parameter values are far from the
correct solution \citep{titter:85,pilla:01}. However, such rules are
widely employed and therefore we conducted an experiment to assess the
two criteria on two data sets. 

The most important assessment of the convergence of an ML algorithm is
the value of the loglikelihood, as it provides information about the
accuracy of parameter estimators on a confidence interval
scale. Therefore, loglikelihood-based criterion is a useful one to
employ in assessing the convergence of an algorithm in finding the MLE
of the parameters \citep{lindsay:95,pilla:01}.

\vspace{0.1in}
\noindent{\sf Simulation Experimental Design:} We consider the
simulated data by generating a sample of size $n = 270$ from $N_p(\mQ,
{\bf I})$ with $p = 3$, where $N_p(\mQ, {\bf I})$ represents a
measure of a $p$-dimensional normal random variable with mean $\mQ$
and an identity variance-covariance matrix. The true mixing measure
for $\mQ \in \mG$, is chosen by selecting the coordinates of
$\boldsymbol{\theta}_j \, (j = 1, \ldots, m)$ from the set $\{-5, 0,
5\}$ in all possible combinations, with equal mass at each support
vector. This resulted in a total of $m = 3^3$ mixture components.

\vspace{0.1in}
\noindent{\sf Fisher Iris Data:} We fit a mixture of multivariate
normal distributions to Fisher iris data \citep{fisher:36}. The data
consists of $n = 150$ observations collected on flowers of three iris
species (Setosa, Verginica and Versicolor). Each observation is a
vector of $p = 4$ variables sepal length (${\bf y}_1$), sepal width
(${\bf y}_2$), petal length (${\bf y}_3$) and petal width (${\bf
y}_4$).

\begin{table}[htp]
\caption{Effect of a convergence criterion on the final
loglikelihood in fitting the fixed support mixture model.}
\label{table:compare0}\vspace{.3cm} 
\centerline{
\begin{tabular}{lrrcc} \hline 
\multicolumn1c{\bf Experiment} & \multicolumn1c{$t$} &
\multicolumn1c{$l\left(\boldsymbol{\pi}^{(t)} \right)$} 
& $\Psi \left(\mQ^{(t)} \right)$ & \multicolumn1c{$\Lambda^{(t)}$} \\ 
 \hline    
Simulated  & 1067 & -2313.6826  &  0.0830 & 0.0536 \\
Fisher Iris & 460 & -376.9595  &  3.0017 & 0.0156 \\ \hline 
\end{tabular} }
\end{table}

For each of the data sets, we selected the observed data matrix ${\bf
y}$ for ${\boldsymbol \Theta}_m$ and also set $\widehat{\boldsymbol
\Sigma} = {\bf S}$, the sample variance-covariance matrix. In order to
assess the accuracy of the algorithms at a given step $t$, we
found the final loglikelihood value $l(\widehat{\boldsymbol{\pi}})$ to
a high degree of accuracy using the PD algorithm for a sufficiently
large $t$. Next, we fit mixtures of multivariate normal distributions
to the simulated and iris data sets via the D-EM algorithm using the
convergence criterion (\ref{eq:lstop}) with $\tau = 0.0001$. The
$\Lambda^{(t)}$ values, presented in Table \ref{table:compare0},
demonstrate that the convergence criterion (\ref{eq:lstop}) would
result in substantially less than four decimal accuracy for the Fisher
iris data. On the other hand, the criterion based on $\Psi(\mQ^{(t)})$
in (\ref{eq:gstop}) guarantees the final accuracy.

\subsection{Empirical Assessment of Convergence Rate}

We empirically assess the rate of convergence of the PD algorithm
relative to the D-EM algorithm by defining $\Lambda^{(t)}$ in
(\ref{delta}) as the {\em residual of the loglikelihood} at the $t$th step.

To be precise, for some $\mbp^{(0)} \in {\boldsymbol \Pi}$, let
$\left\{\mbp^{(t)} \right\}_{t \geq 1}$ be a sequence in ${\boldsymbol
\Pi}$ generated by an algorithm (such as the PD and D-EM
algorithms). The algorithm can be expressed as $\mbp^{(t)} \in \mM
\left(\mbp^{(t - 1)} \right)$ for $t \geq 1$, where the map $\mM:
{\boldsymbol \Pi} \to 2^{{\boldsymbol \Pi}}$ is a point-to-set
mapping. If $\mbp^{(t)}$ converges to $\wmbp$ and $\mM(\cdot)$ is
continuous, then $\wmbp$ must satisfy $\wmbp \in
\mM(\wmbp)$.

\vspace{0.1in}
\noindent{\sf Definition 2 (Asymptotic Convergence Rate):} Assume
that $\wmbp = \mM(\wmbp)$ and that the sequence $\{\mbp^{(t)}\}_{t
\geq 1}$ is generated by the map $\mM$ such that $\lim_{t \to \infty}
\, \mbp^{(t)} = \wmbp$. Under the regularity conditions given by
\cite{wu:83}, this implies that $\lim_{t \to \infty} \, l
\left(\mbp^{(t)} \right) = l(\wmbp)$. The {\em asymptotic convergence
rate of the loglikelihood sequence} $\left\{l \left(\mbp^{(t)} \right)
\right\}_{t \geq 1}$ at $l(\wmbp)$ generated by an algorithm is
defined as 
\bes
  r := \lim_{t \to \infty} \; \big| l \left(\wmbp \right) - l
  \left(\mbp^{(t)} \right) \big|^{\frac{1}{t}} 
\ees 

From the following lemma \citep{pilla:01}, the smaller the $r$ for any
given loglikelihood sequence, the faster it is progressing towards the
MLE. 
\begin{lemma}
  \label{lem-rate} {\rm If the sequence $\{l(\mbp^{(t)})\}_{t \geq 1}$
  is converging linearly, then as $t \to \infty$, the slope of the
  curve obtained by plotting $\log \, \{\Lambda^{(t)}\}$ against $t$
  converges to $\log \, (r)$, where $r$ is the asymptotic rate of
  convergence of the loglikelihood sequence generated by an
  algorithm.}
\end{lemma} 

\begin{figure}[ht] 
  \vspace{-0.1in}
  \scalebox{0.5}{\includegraphics{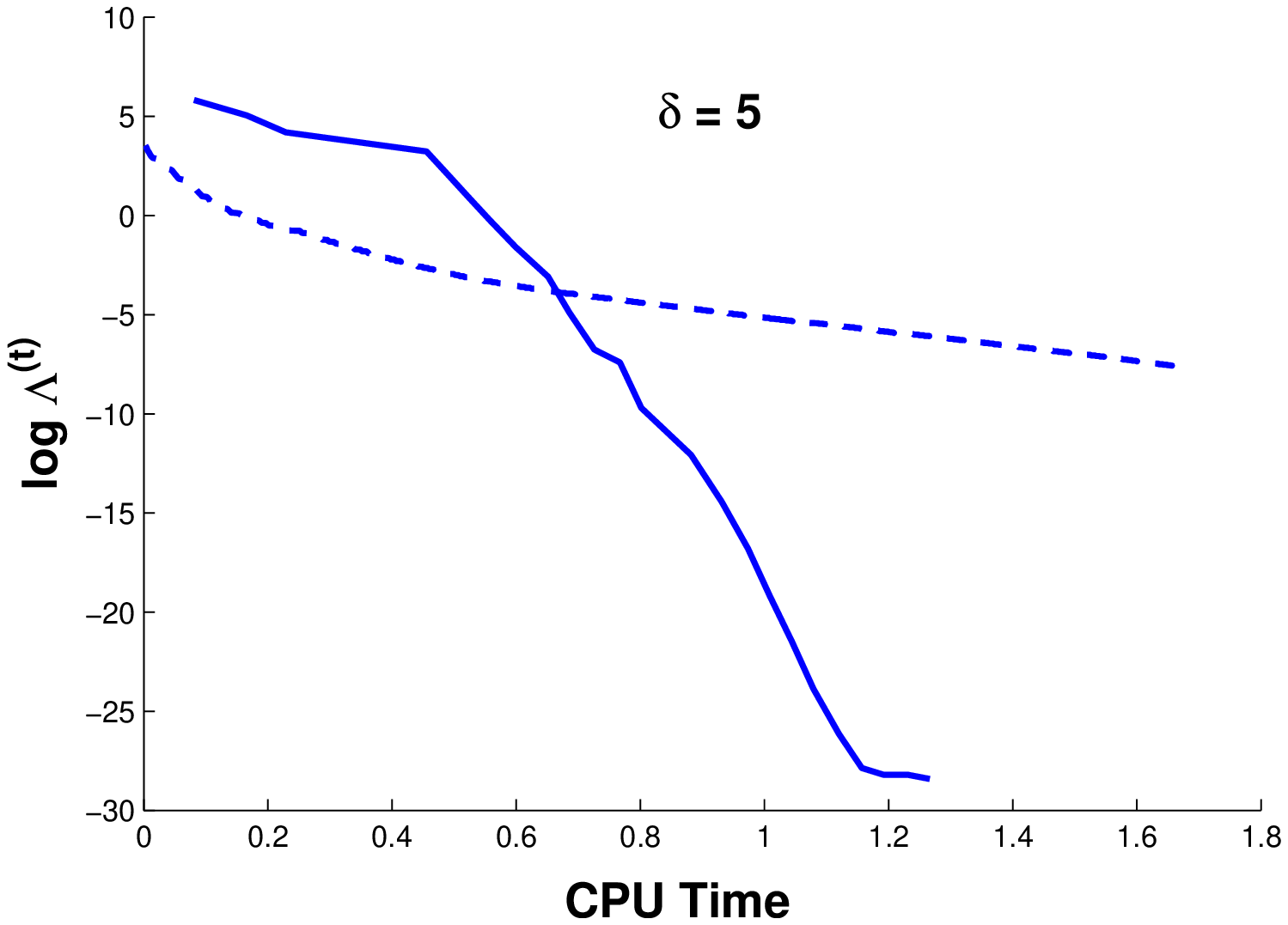}} 
  \scalebox{0.5}{\includegraphics{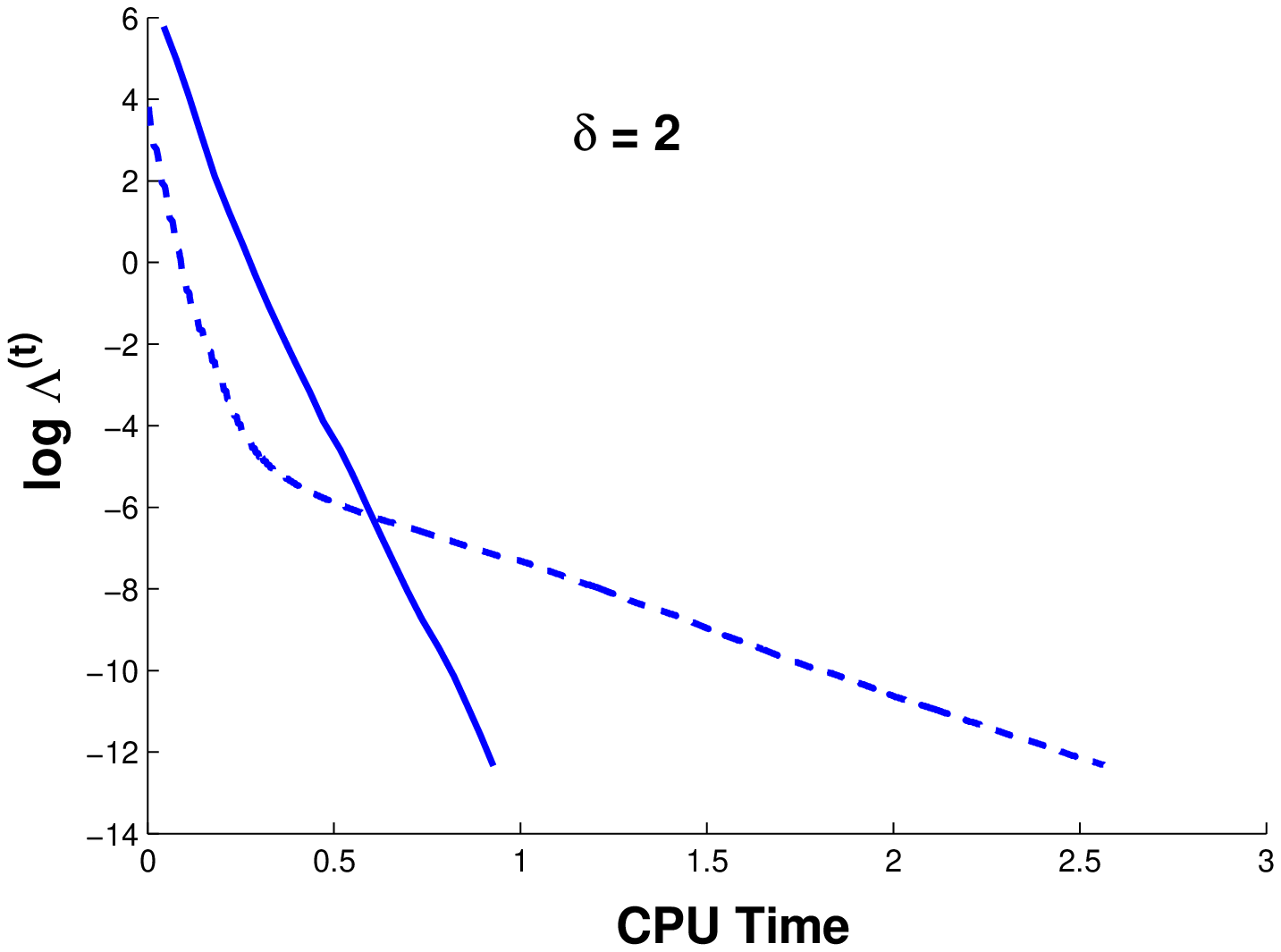}}  \\ 
  \scalebox{0.5}{\includegraphics{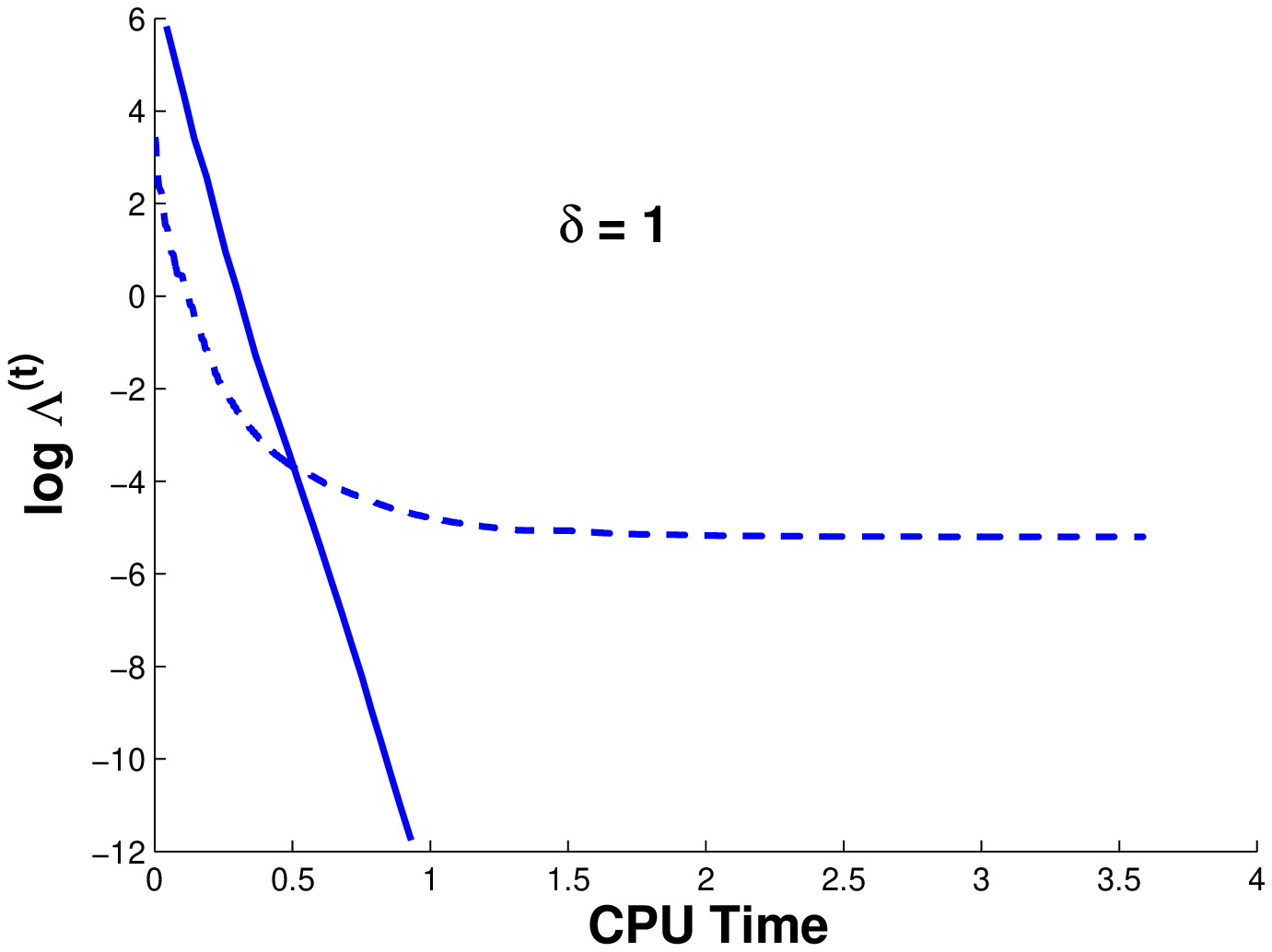}} 
  \scalebox{0.5}{\includegraphics{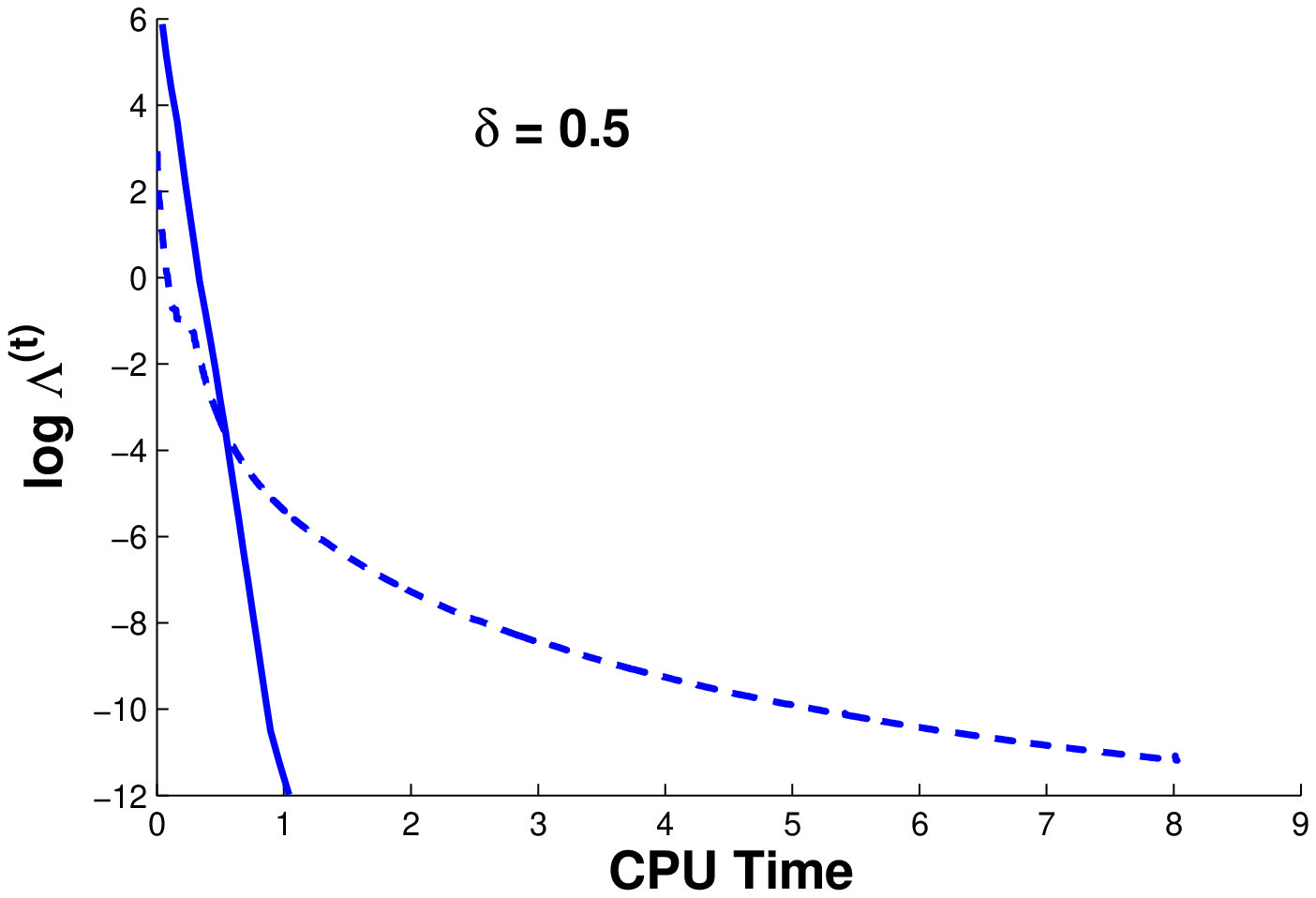}} 
  \vspace{-0.05in}
\caption{Plot of the log residual of the loglikelihood, $\log \,
\{\Lambda^{(t)}\}$ against CPU time for the Fisher iris data with four
variables demonstrating the sublinear convergence of the discrete EM
(dashed line) and the linear convergence of the Penalized Dual (solid
line) algorithms.}
\label{fig-rate}
\end{figure}

In order to assess the rate of convergence of the PD, relative to the
D-EM, algorithm, we consider the Fisher iris data considered earlier
for fitting a collection of semiparametric mixture of multivariate
normal distributions $N_p(\mQ, \, \delta \, {\boldsymbol \Sigma})$,
where $p = 4, {\boldsymbol \Sigma} \in \mS$ and $\delta \in \Re_+$
(details in Section \ref{sec-tree}). As before, we selected the
observed data ${\bf y}$ for ${\boldsymbol \Theta}_m$ and set $\wS =
{\bf S}$ in fitting the PD and D-EM algorithms. Figure
\ref{fig-rate} demonstrates the behavior of the algorithms for $\delta
\in \{5, 2, 1, 0.5\}$. We used logarithmic scaling of the vertical
axis since a linearly convergent algorithm will become linear on this
scale as $t \to \infty$. Note that for the fixed support mixture
model, the D-EM algorithm is converging sublinearly whereas the PD is
converging linearly to the MLE $\wmbp$; a significant improvement in
convergence rate. In fact, \cite{pilla:01} observed a similar behavior
of sublinear convergence of the D-EM algorithm for a class of
univariate finite mixture problems.

\section{Semiparametric Mixtures of Multivariate Normal
Distributions} 
\label{sec-struc}

The methodology developed in this article is applicable to a wide
range of problems, including multivariate {\em t} mixtures. However,
the particular interest here is in difficult problems with
multivariate normal mixtures due to its ubiquitous applications. In
semiparametric mixture setting, the mixing distribution $\mQ$ is
modeled nonparametrically in the presence of an unknown ${\boldsymbol
\Sigma} \in \mS$, the variance-covariance matrix common to all $m$
components. 
 
\subsection{Structural Properties} 
\label{sec-norm}

Let ${\rm g}_{_\mQ}({\bf y}_i; {\boldsymbol \Sigma}) = \sum_{j = 1}^m
\pi_j \, f_{\boldsymbol{\mu}_j}({\bf y}_i; {\boldsymbol \Sigma})$ for
${\bf y}_i \in \mY$ be a finite mixture of $p$-dimensional normal
distributions, where $\boldsymbol{\mu}_j \in \Re^p \subset \Omega$ is
the mean vector of the $j$th component density
$f_{\boldsymbol{\mu}_j}({\bf y}_i; {\boldsymbol \Sigma})$ and 
${\boldsymbol \Sigma} \in \mS$ is common to all $m$ components. Note
that in the continuous case $d = n$. The corresponding loglikelihood
is expressed as
\bea
  \label{eq:lQS} 
  l(\mQ; {\boldsymbol \Sigma}) = \sum_{i = 1}^n \ln \, \{{\rm
  g}_{_\mQ}({\bf y}_i; {\boldsymbol \Sigma})\}. 
\eea

In the univariate case, \cite{char:05} establish that for a general
family of mixture models with a structural parameter $\beta$ (e.g.,
$\sigma^2$ in the normal case), the likelihood framework breaks down
when joint estimation of $m, \mQ$ and $\beta$ is attempted: at best
the joint estimator of $m, \mQ$ and $\beta$ is degenerate, and at
worst it does not even exist. The ML fails in this setting since
taking finite samples from continuous probability distributions yields
discrete data sets.  When models that closely mimic discrete
probability distributions are available, as they are when there are no
restrictions on $\mQ$ and $\beta$, the likelihood will favor such
models. The NPMLE results of Lindsay (1995, Section 2.6) cannot be
applied if ${\boldsymbol \Sigma}$ is unknown; however, the following
result holds.

We define the gradient function for the multivariate normal mixture
distributions as
\begin{eqnarray}  
  \label{eq:ngrad}
  \mD_{\widehat{\mQ}_{_{\boldsymbol \Sigma}}}(\boldsymbol{\mu};
  {\boldsymbol \Sigma}) := \sum_{i = 1}^n \left\{
  \frac{f_{\boldsymbol{\mu}}({\bf y}_i; {\boldsymbol
  \Sigma})}{{\rm g}_{_{\wQ_{\boldsymbol \Sigma}}}({\bf y}_i; 
  {\boldsymbol \Sigma})} - 1 \right\} \quad \mbox{for}
  \quad \boldsymbol{\mu} \in \Omega.
\end{eqnarray}
Next, we define an NPMLE of $\mQ \in \mG$ for a fixed $\bS \in \mS$ as
\bes
  \wQ_{_{\boldsymbol \Sigma}} = \underset{\mQ \, \in \, \mG}{\arg \,
  \max} \; l(\mQ; \bS).  
\ees

\begin{theorem}[Unique NPMLE of $\mQ$]  
\label{thm-mmml}
{\rm Assume ${\boldsymbol \Sigma} > 0$ is fixed. \\ (1) Suppose
$\wQ_{_{\boldsymbol \Sigma}}$ satisfies
\bea
  \label{gradineq}
  \mD_{\widehat{\mQ}_{_{\boldsymbol \Sigma}}}(\boldsymbol{\mu};
  {\boldsymbol \Sigma}) \leq 0 \quad \mbox{for all} \quad
  \boldsymbol{\mu} \in \Omega, 
\eea
  then $\wQ_{_{\boldsymbol \Sigma}}$ is an NPMLE of $\mQ \in
  \mG$. \\ (2) Let the set $\{\boldsymbol{\mu}_1, \ldots,
  \boldsymbol{\mu}_K\}$ for some $K \leq n$ be the solution set
\bes
  \label{cond1}
  \left\{ \boldsymbol{\mu}: \mD_{\widehat{\mQ}_{_{\boldsymbol
  \Sigma}}}(\boldsymbol{\mu}; {\boldsymbol \Sigma}) = 0 \right\}.
\ees
If the vectors
\bes
  \label{eq:fvec} 
  {\bf  f}_{\boldsymbol{\mu}_j}({\bf y}; {\boldsymbol \Sigma}) =
  \left\{f_{\boldsymbol{\mu}_j}({\bf y}_1; {\boldsymbol \Sigma}),
  \ldots, f_{\boldsymbol{\mu}_j}({\bf y}_n; {\boldsymbol \Sigma})
  \right\}^T \; \mbox{for} \; j = 1, \ldots, K 
\ees
are linearly independent, then $\widehat{\mQ}_{_{\boldsymbol \Sigma}}$
is the unique NPMLE of $\mQ \in \mG$.}
\end{theorem} 

For the multivariate normal mixture model with a common ${\boldsymbol
\Sigma}$, we restrict attention to finite discrete latent
distributions $\mQ$, then the pair $(\mQ, {\boldsymbol \Sigma})$ is
identifiable \citep{lindsay:95}. For a general family of univariate
mixtures, \cite{char:05} establish that joint estimation of $m, \mQ$
and $\beta$ is a {\em well-defined} problem if $\mQ$ is finitely
supported. If $\mQ$ is not finitely supported, then ${\rm
g}_{_\mQ}({\bf y}; {\boldsymbol \Sigma})$ need not determine $\mQ$ and
$\beta$ uniquely. Hence, an ML approach to the joint estimation of $m,
\mQ$ and ${\boldsymbol \Sigma}$ fails. However, since we fix $m$ and
consider $\mQ$ to be finitely supported, joint estimation of $\mQ$ and
${\boldsymbol \Sigma}$ is feasible. Therefore, we can apply Algorithm
1 described in Section \ref{sec-parmix} to jointly estimate $\mQ$ and
${\boldsymbol \Sigma}$. 

The following theorem establishes that joint identifiability of
$(\mQ, {\boldsymbol \Sigma})$ fails if $\mQ$ is not finitely
supported. The proof follows from the univariate nesting structure
result, under mild regularity conditions, given by \cite{char:05}.

\begin{theorem}[Multivariate Mixture Nesting Structure] 
\label{thm-nest} 
{\rm The class of multivariate normal mixture distributions possesses
the nesting structure. That is, for any ${\boldsymbol \Sigma}^{\dag}
\succ {\boldsymbol \Sigma}$, in the sense of L\"{o}wner ordering, 
\begin{eqnarray*}
   \{N_p(\mQ, {\boldsymbol \Sigma}^{\dag})\!: \mQ \in \mathcal{G}\}
   \subseteq \{N_p(\mQ, {\boldsymbol \Sigma})\!: \mQ \in
   \mathcal{G}\},   
\end{eqnarray*}
where $N_p(\mQ, {\boldsymbol \Sigma})$ represents a measure of a
$p$-dimensional normal random variable with mean $\mQ$ and a
variance-covariance matrix ${\boldsymbol \Sigma} \in \mS$. }
\end{theorem} 

\subsection{Role of the Sieve Parameter in Building Semiparametric
Mixture Models}     
\label{sec-tree}

We investigate building the sieve of models $N_p({\bf F}, \, \delta \,
{\boldsymbol \Sigma})$, where ${\boldsymbol \Sigma} \in \mS$ and
$\delta \in \Re_+$ is a {\em sieve parameter} (similar to the
smoothing parameter employed in density estimation). The sieve
parameter controls the dimensionality of a mixture model as will be
demonstrated later. We derive theory for building a collection of
semiparametric mixture models, including the multivariate case.

\begin{figure}[htp] 
  \centerline{\scalebox{0.5}{\includegraphics{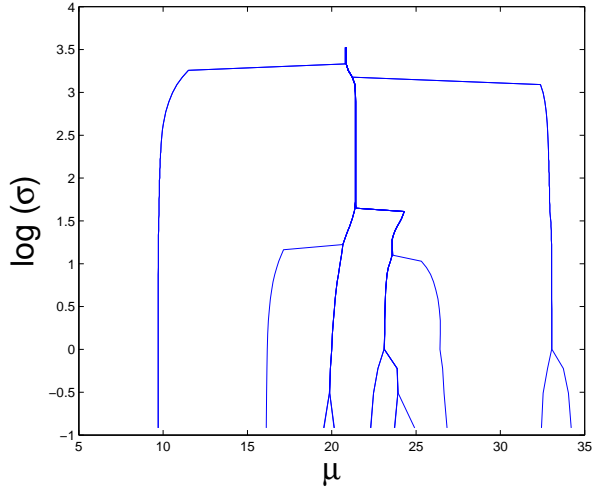}} }
\caption{A Mixture Tree for the galaxy data set.}
\label{fig-astrod}
\end{figure}

\begin{figure}[htp] 
  \centerline{\scalebox{0.75}{\includegraphics{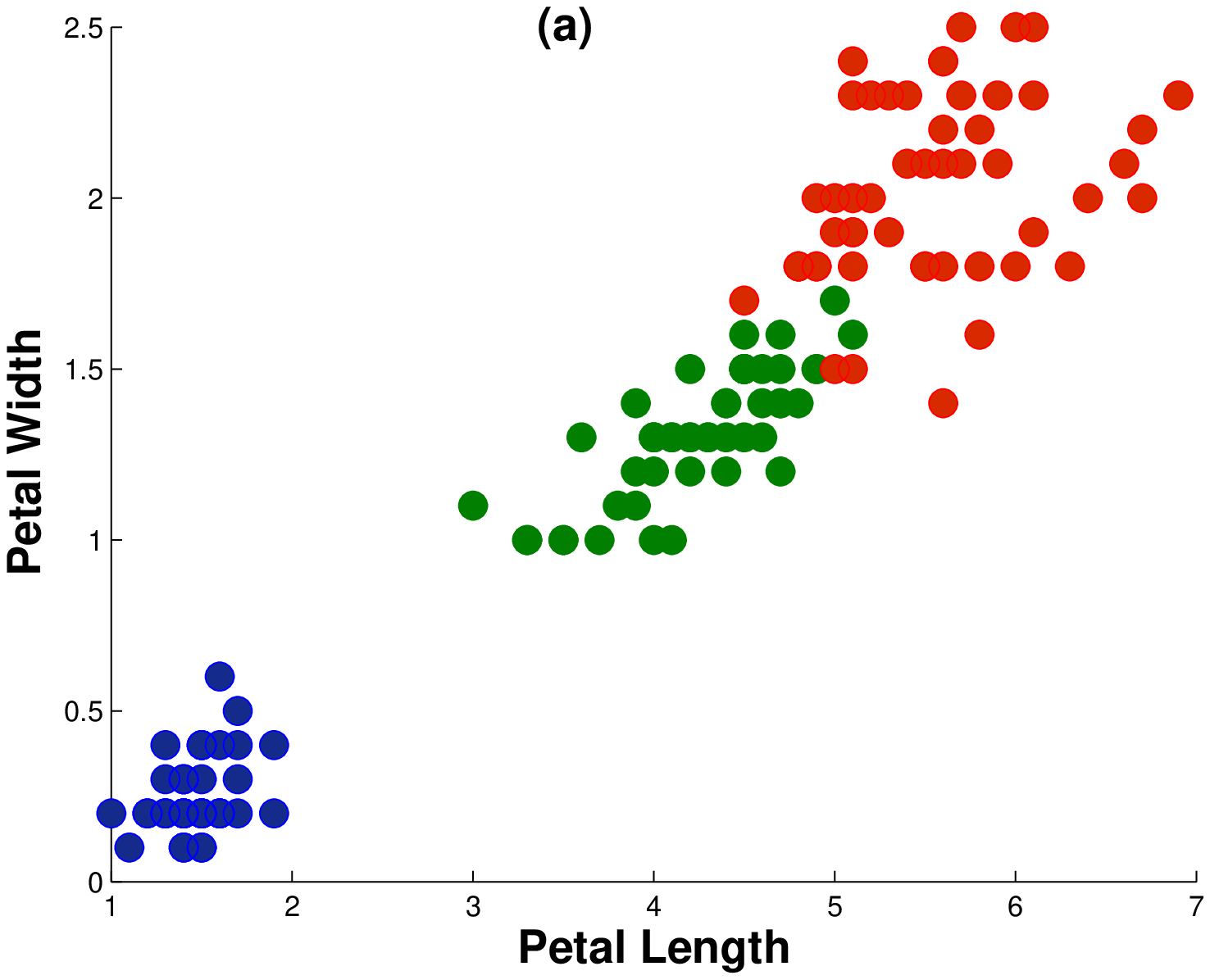} } } 
\vspace{-0.1in}
  \centerline{\scalebox{0.75}{\includegraphics{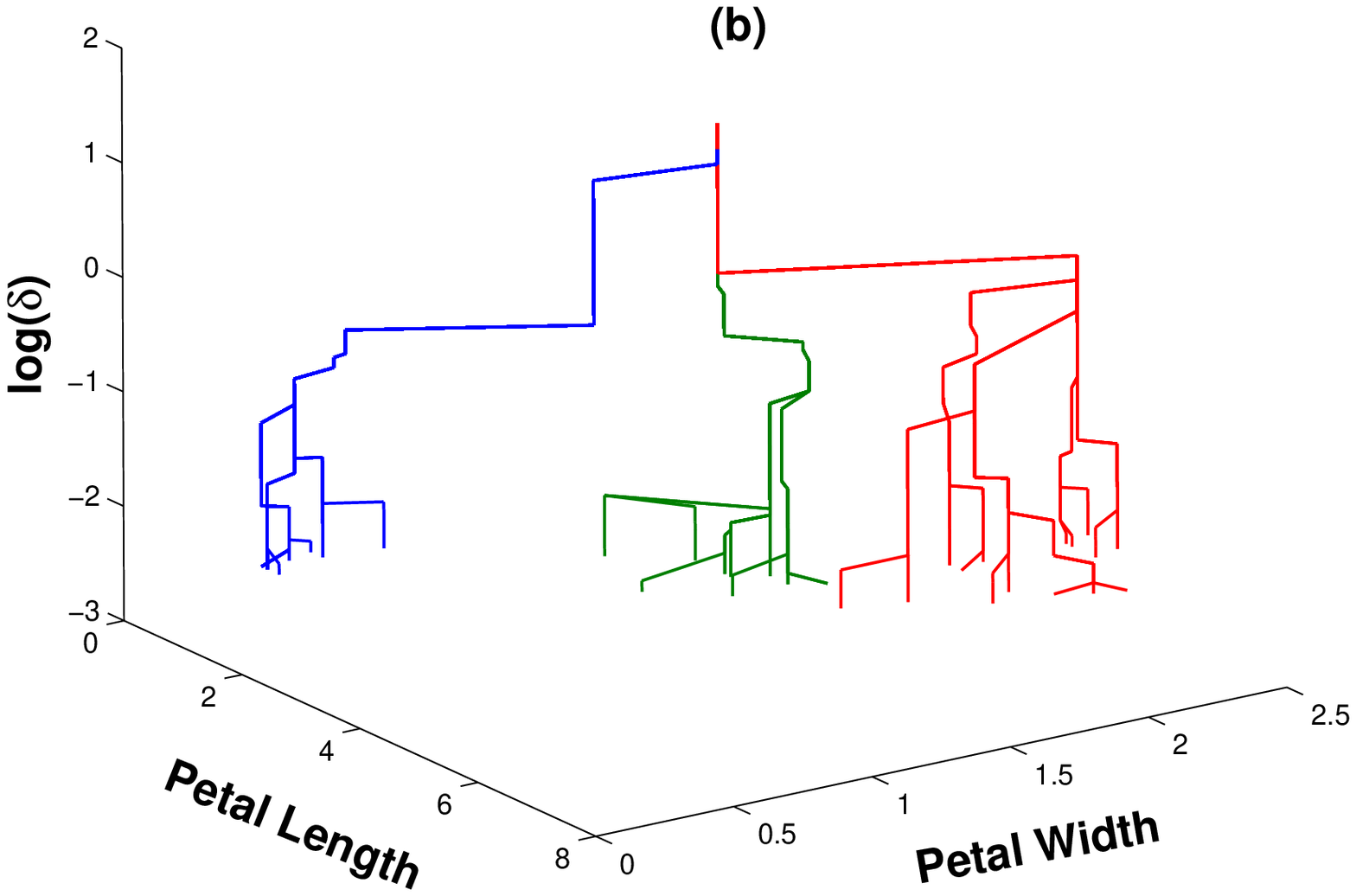} } }
\caption{(a) The scatter plot of the Fisher iris data for the two
variables, namely the petal length and petal width, showing three main
groups. (b) A Mixture Tree demonstrating the effect of the sieve
parameter $\delta$ on the number of components. The three main
branches in the tree correspond to the three main components in (a).}
\label{fig-irisd}
\end{figure}

In order to create a general family of mixture models, we consider the
class $\{N_p(\mQ, \, \delta \, {\boldsymbol \Sigma}): \mQ \in
\mathcal{G}, {\boldsymbol \Sigma} \in \mS, \delta \in
\Re_+\}$. For $\delta_1 > \delta_0$, Theorem \ref{thm-nest} implies
that
\begin{eqnarray*}
   \{N_p(\mQ, \, \delta_1 \, {\boldsymbol \Sigma})\!: \mQ \in
   \mathcal{G}\} \subseteq \{N_p(\mQ, \, \delta_0 \, {\boldsymbol
   \Sigma})\!: \mQ \in \mathcal{G}\};
\end{eqnarray*}
hence, the collection of models becomes richer as $\delta \rightarrow
0$ (see Figures \ref{fig-astrod} and \ref{fig-irisd}). Moreover, every
$p$-dimensional distribution ${\bf F}$ can be obtained as the weak
limit of $N_p({\bf F}, \, \delta \, {\boldsymbol
\Sigma})$ as $\delta \rightarrow 0$. Therefore,
we can approximate any distribution by choosing $\delta$ small. As a
consequence, the principle of maximum likelihood cannot be applied to
select $\delta$ in the model $N_p(\mQ, \, \delta \, {\boldsymbol
\Sigma})$ since the likelihood becomes unbounded as $\delta
\rightarrow 0$. \cite{char:05} develop theory for the univariate
mixtures and demonstrate the effect of small $\delta$ for a general
family of univariate mixtures which extends to the multivariate case
considered here.

We create a strategy for building a collection of models
$\{(\mQ_{\delta}, \, \delta \, {\boldsymbol \Sigma})\!: \mQ 
\in \mathcal{G}, \delta \in \Re_{+} \}$ using the Penalized
Dual algorithm. As $\delta \rightarrow 0$, the NPMLE
$\widehat{\mQ}_{\delta}$ converges in distribution to $n^{-1} \sum_i
\vartheta({\bf y}_i)$, where $\vartheta({\bf y})$ is a discrete
measure concentrated at ${\bf y}$.

To demonstrate the effect of $\delta$ on the mixture complexity, we
create a collection of models for both the univariate and multivariate
data. In the univariate case, we simply have a $\sigma$ parameter. The
univariate application considers the galaxy data set [Table 1 of
\cite{roeder:90}] of 82 observations of relative velocities for
galaxies from six well separated conic sections of the Corona Borealis
region. Scientific interest lies in identifying substructures in
clusters of galaxies. Multimodality is evidence of voids and
superclusters in the far universe. \cite{roeder:90} obtained
$\widehat{\sigma} = 0.95$ using least squares cross validation. We set
$\mbm \in {\boldsymbol \Theta}_m = \{9, \ldots, 35\}$ with a grid size
of 0.02 for building a collection of semiparametric mixture models
using Algorithm 2. The plot of $\log \, (\sigma)$ against the support
set $\mbm$ corresponding to the estimate $\wmbp$ obtained using the PD
algorithm (namely, Algorithm 2) is shown in Figure
\ref{fig-astrod}. That is, at each fixed $\log \, (\sigma)$, the plot 
displays $\mu$ parameter values that have positive mixture
probability. The figure demonstrates the effect of $\sigma$ on the
mixture complexity $m$.

Next, we consider the Fisher's iris data described earlier. Once
again, we selected observed data ${\bf y}$ for ${\boldsymbol
\Theta}_m$ and set $\wS = {\bf S}$. We let $\delta = \{0.1, \ldots,
5\}$ for building a collection of semiparametric multivariate mixture
models using Algorithm 2. Figure \ref{fig-irisd} shows the effect of
$\delta$ on the mixture model complexity when only the two variables,
namely the petal length and petal width are considered. The galaxy and
Fisher iris data sets demonstrate that the number of components is a
consequence of the choice of $\sigma$ (or $\delta$ as the case may be)
rather than a pre-selected parameter.

\vspace{-0.1in}

\subsection{Selection of the Support Set of $\mQ$}
\label{sec-grid}

As described in Section \ref{sec-parmix}, in approximating
$\mathcal{G}$, the biggest challenge is in selecting a suitable
${\boldsymbol \Theta}_{m}$ while keeping computations manageable. This
is addressed in this section.

In the absence of a prior knowledge of the mixture complexity, correct
specification of ${\boldsymbol \Theta}_m \subset \Omega$, the support
set of $\mQ$, is very important for the Step 1 of Algorithm 1 (or
equivalently for Algorithm 2). As expected, the final loglikelihood
depends on this choice. In this section, we illustrate through the
simulated data described earlier how the observed data matrix ${\bf
y}$ provides the best choice for approximating the continuous
parameter space $\Omega$. In effect, we select
$\{\boldsymbol{\theta}_1 = {\bf y}_1, \boldsymbol{\theta}_2 = {\bf
y}_2, \ldots, \boldsymbol{\theta}_m = {\bf y}_d\}$. Note that choosing
${\bf y}$ for the discrete parameter space ${\boldsymbol \Theta}_m$
clearly covers the region of likely support vectors for the normal
means and has the advantage of adapting naturally in richness to the
sample size of the problem.

To assess the effectiveness of using ${\bf y}$ for ${\boldsymbol
\Theta}_m$ (which is approximating the continuous parameter space
$\Omega$) we consider the simulated data described in Section
\ref{sec-criteria}.  The true mixing measure for $\mQ$ chosen for the
simulation experiment is denoted by ``True Support'' in Table
\ref{table:compare}. The ``Equi-Distant'' set for ${\boldsymbol
\Theta}_m$ was constructed on a lattice by choosing the elements in
$\boldsymbol{\theta}_j = (\theta_{j1}, \theta_{j2},
\theta_{j3})$ for $j = 1, \ldots, 8^3$ from the set $\{-7, -5, \ldots,
5, 7 \}$ resulting in a total of $m = 8^3$ support vectors; this set
also included all the true support vectors. Table
\ref{table:compare} presents results obtained using Algorithm 2
(i.e., fixed support mixture model of estimating $\mbp$ for a given
${\boldsymbol \Theta}_m$) and Step 2 of Algorithm 1 (i.e., continuous
support mixture model of estimating $\mQ$ and ${\boldsymbol \Sigma}$
for a fixed $m$). 

In general, ${\bf y}$ should be an effective choice for ${\boldsymbol
\Theta}_m$ given that equi-distant is still a subjective one in the
absence of any knowledge about the length of the distance. From the
theory presented in Section \ref{sec-parmix}, as $m \to \infty$,
${\boldsymbol \Theta}_m \to \Omega$. However, in practice, choosing $m
= n$ is effectively creating a dense set for ${\boldsymbol \Theta}_m$
and in fact approximating $\Omega$ very well.

\begin{table}[htp]
\caption{Effect of $\boldsymbol{\Theta}_m \subset \Omega$, the support
set for $\mQ$, on the estimated loglikelihood. We set
$\widehat{\boldsymbol \Sigma} = {\bf S}$ in finding $l(\wmbp)$ using
the Penalized Dual (PD) algorithm. The solution $\wmbp$ obtained from
the PD algorithm with the corresponding fixed support set and ${\bf
S}$ are employed as parameter starting values for finding
$l\left(\wQ_{\delta}, \delta \, \wS \right)$ for a fixed $\delta =
0.2$ using the continuous EM (C-EM) algorithm.} 
\label{table:compare} 
\vspace{.3cm} 
\centerline{
\begin{tabular}{lcc} \hline 
\multicolumn1c{${\boldsymbol \Theta}_m$} &
$l\left(\widehat{\boldsymbol{\pi}} \right)$ & $l \left(\wQ_{\delta},
\delta \, \wS \right)$ \\   
\hline True Support & -2181.9 & -1936.9 \\
Equi-Distant & -2182.8 & -1901.7 \\
Observed Data & -2178.6 & -1876.0 \\ \hline
\end{tabular} }
\end{table}

\section{Applications and Simulation Experiment}
\label{sec-appl}

The applications in this section are used to investigate the roles of
many overlapping components which create an ideal situation for solving
the large-scale practical problems. 

We assess the performance of the algorithms in finding the NPMLE of
$\mQ$ and for fitting the collection of semiparametric mixture models
with applications to several data sets. The data sets, the parameter
estimates and the {\sf Matlab} software for fitting mixtures are
available from the first author. For the Step 1 of Algorithm 1, we 
set $\widehat{\boldsymbol \Sigma} = {\bf S}$; however, we estimate it
in Step 2.

\subsection{Mortality Data}
\label{sec-poiss}

Our first application illustrates the tremendous advantage of our
method in reducing the dimension of discrete mixture problems. In
these problems the magnitude of $d$, the number of distinct observed
data points, could be much smaller than $m$, the cardinality of
${\boldsymbol \Theta}_{m}$.

We consider the data on death rates which gives the number of death
notices for women aged 80 and over, from the \textit{Times} newspaper
for each day in the three-year period 1910 to 1912
\citep{titter:85}. For the later data sets, we chose the observed 
data matrix ${\bf y}$ as the support set ${\boldsymbol
\Theta}_m$. However, for this application, we selected the support set
to be ${\boldsymbol \Theta}_m = \{ 0, 0 + \eta, \ldots, 9 - \eta, 9
\}$, where $\eta \in \{1, 0.5, 0.1, 0.01\}$. In effect, the mixture
complexity $m \in \{10, 20, 100, 1000\}$. It is worth noting that the
dimension of the dual optimization problem is $d$ (equals 10) whereas
that of the mixture problem is $(m - 1)$ which grows significantly
with the cardinality of the set $\eta$.

\begin{table}[htp]
\small{
\caption{Building a collection of finite mixture models
$\left\{\wQ_{\eta}\!: \eta \in \Re_{+} \right\}$ using Algorithm 1 for
the mortality data. First step involved setting ${\boldsymbol
\Theta}_m = \{0, 0 + \eta, \ldots, 9 - \eta, 9 \}$ with $\eta \in \{1,
0.5, 0.1, 0.01\}$ and finding $l(\wmbp)$ using the PD-based and
discrete EM (D-EM) algorithms. The estimate $\wmbp$ obtained from the
PD algorithm with the corresponding fixed support set is used as
parameter starting values for finding $l\left(\wQ_{\eta} \right)$
using the continuous EM (C-EM) algorithm.}
\label{table:compare4} \vspace{.3cm} 
\centerline{
\begin{tabular}{lccccrr} \hline 
& & \multicolumn1c{$l(\widehat{\boldsymbol{\pi}})$} &
\multicolumn1c{$\Lambda^{(t)}$} & \multicolumn1c{$\Psi \left(\mQ^{(t)}
\right)$} & & \multicolumn1c{\bf CPU} \\ {\bf Algorithm} & $\eta$ &
\multicolumn1c{$l \left(\wQ_{\eta} \right)$} & \multicolumn1c{$\times
10^{3}$} & \multicolumn1c{$\times 10^{3}$} & \multicolumn1c{N$^{(t)}$} &
\multicolumn1c{\bf Factor} \\  
\hline 
PD & 1 &   -1990.0928  &   0.0000  &   0.2577  &   25 &   5 \\ 
PD$^{\rm IC}$ & &   -1990.0928  &   0.0000  &   0.2577  & 25 &
7 \\ 
D-EM & &   -1990.0929  &   0.0172  &   4.9885  &   1,238 &  1 \\ 
C-EM & & -1989.9 & - & - & 2,179 & - \\ [2ex] 
PD & 0.5 &   -1989.9941  &   0.0000  &   0.1881  &   26 &
120 \\ 
PD$^{\rm IC}$ & &   -1989.9941  &   0.0000  &   0.1881  &   26  &
142 \\ 
D-EM &  &   -1989.9949  &   0.7136  &   4.9997  & 31,149  &   1 \\  
C-EM & & -1989.9 & - & - & 2,360 &- \\ [2ex] 
PD & 0.1  &   -1989.9281  &   0.0000  &   0.2521  &   25 &
638 \\ 
PD$^{\rm IC}$ & &   -1989.9281  &   0.0000  &   0.2520  &   25 &
719 \\ 
EM &  &   -1989.9322  &   4.0901  &   5.0000  &   108,312 &   1 \\  
C-EM & & -1989.9 & - & - & 1,997 & - \\ [2ex] 
PD & 0.01 &   -1989.9272  &   0.1108  &   0.2270  &   27  & 943 \\  
PD$^{\rm IC}$ &  &   -1989.9272  &   0.1108  &   0.2269  &   27  &
1,192 \\ 
D-EM  &  &  -1989.9319  &   4.8230  &   5.0000  &   113,081 &   1 \\   
C-EM & & -1989.9 & - & - & 1,924 & - \\ [2ex] 
\hline 
\end{tabular} } }
\end{table}

We fit a mixture of Poisson distributions to the mortality data using
the PD and D-EM algorithms. Table \ref{table:compare4} presents
N$^{(t)}$, the number of steps required for convergence [based on the
criterion (\ref{eq:gstop})] and ``CPU Factor'', the ratio of the CPU
time required by the D-EM algorithm to that of the PD. This ratio
indicates the factor by which the D-EM algorithm is accelerated. We
also present the values of $l(\wmbp), l(\wQ_{\eta}),
\Lambda^{(t)}$ and $\Psi(\mQ^{(t)})$. Furthermore, we consider the
effect of eliminating the inactive constraints in the PD algorithm,
namely PD$^{\mathrm{IC}}$. The table demonstrates that the PD-based
algorithms advance toward the maximum more rapidly than does the D-EM
algorithm with gains increasing as $m$ (equivalently, the number of
parameters to estimate) increases. Thousand-fold improvements are
obtained at $\eta = 0.01$ for which the number of parameters to
estimate is the largest. For comparison, we fit the same model with
the Rotated EM (an accelerated version of the EM applicable only for
univariate mixtures) developed by \cite{pilla:01} and obtained CPU
factors for the PD, relative to the Rotated EM, as 1.4, 22, 43 and 28,
respectively for $\eta \in \{1, 0.5, 0.1, 0.01\}$.

\begin{figure}[htp] 
  \psfrag{Cumulative Q(theta)}[l][Bl]{{\bf CDF of} $\widehat{\mQ}$} 
  \psfrag{theta}[l][Bl]{$\mbt$} 
  \scalebox{1}{\includegraphics{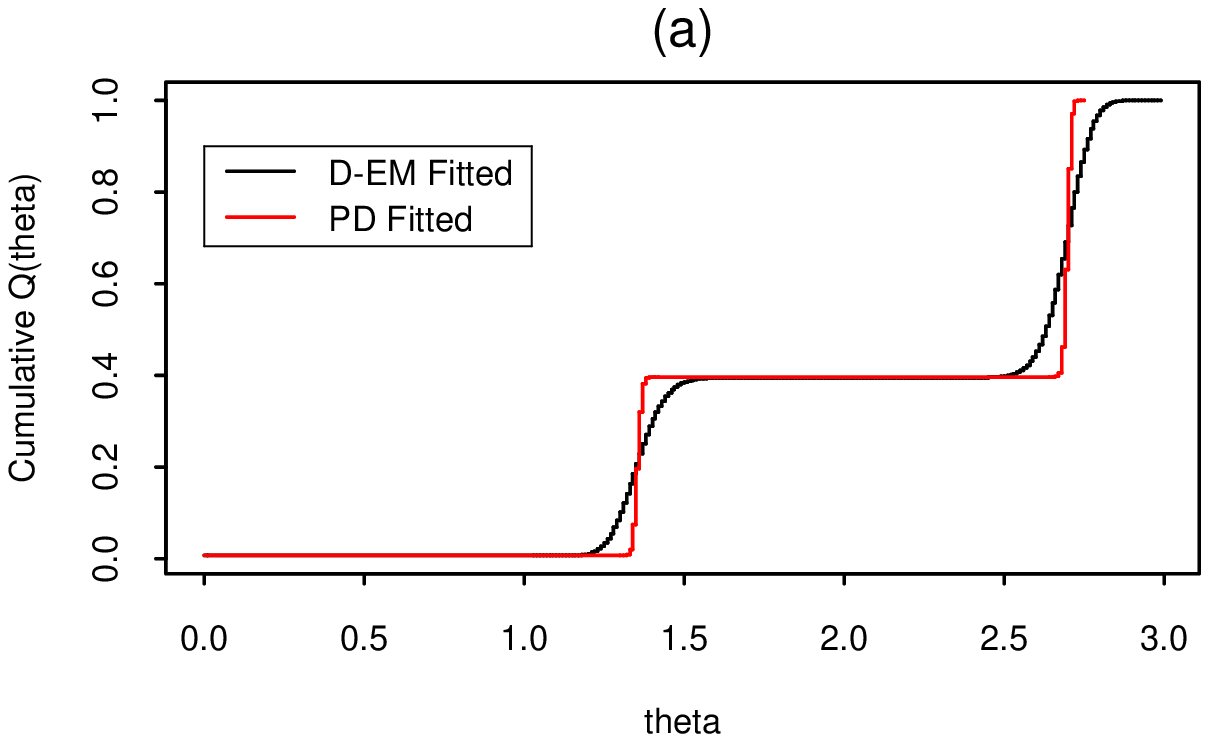}} \\ \\
  \scalebox{1}{\includegraphics{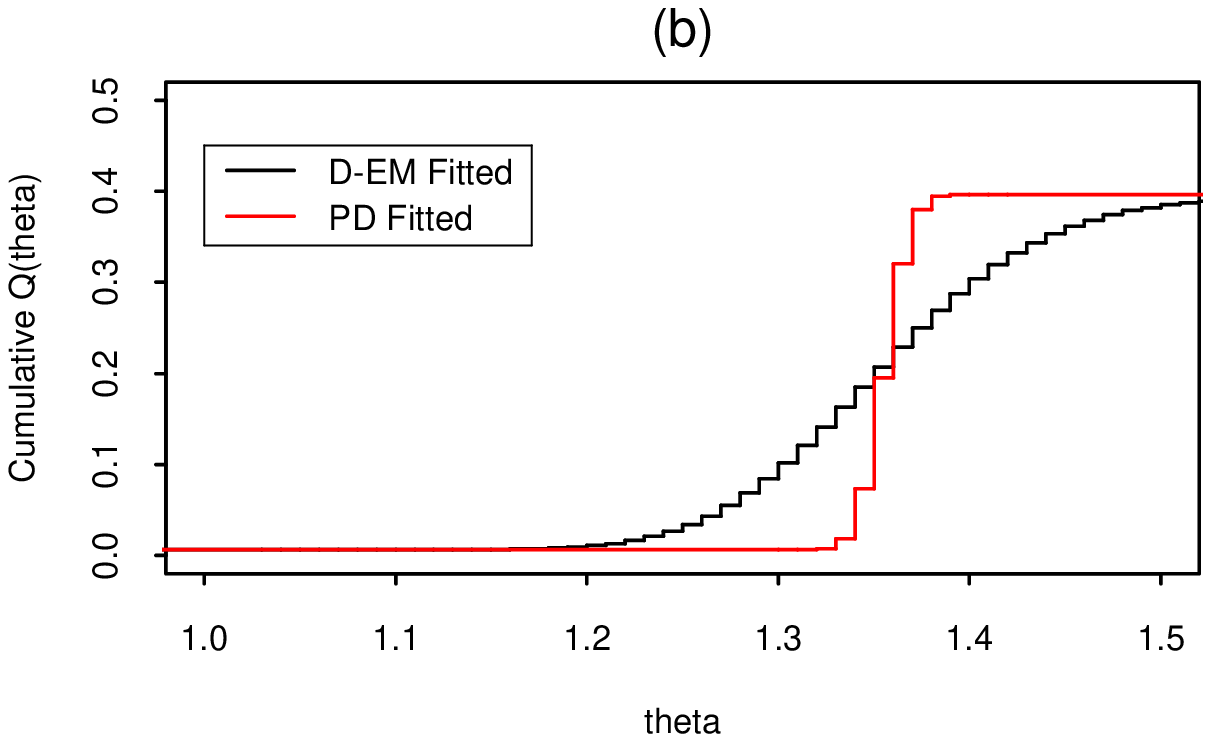}}
\vspace{-0.1in}
\caption{(a) Behavior of the cumulative distribution function (CDF) of
$\widehat{\mQ}$ for the fixed support mixture model obtained from the
D-EM (black) and PD (red) algorithms at $\eta = 0.01$ for the
mortality data. (b) An enlarged view of graph (a) at the first jump;
the D-EM algorithm has many small jumps and retained the redundant 
components.}
\label{fig1}
\end{figure}

At $\eta = 0.01$, the D-EM has retained 199 support points with non
zero probability at convergence (obtaining a smaller loglikelihood
value) whereas the PD has retained just 26 support points and reached
the MLE in a reasonable number of steps; a significant reduction in
the mixture complexity.  For this case of $\delta$, most of the
mixture probabilities are near zero; hence the algorithms must push
the estimates to the boundary of the parameter space---a least
favorable case for the D-EM algorithm. When the NPMLE has fewer than
$d$ support points (an overparameterized mixture problem), then the
D-EM algorithm has great difficulty in allocating probability to the
redundant support points. The behavior of the cumulative distribution
function (CDF) of $\widehat{\mQ}$ for the two algorithms at $\eta =
0.01$ is shown in Figure \ref{fig1}. It is clear that the D-EM
algorithm has an extremely small step size whereas the PD has a
reasonably large step size. This is due to the fact that the D-EM has
retained a significantly large number of components with small
jumps---an artifact of its failure to converge to the MLE in finite
number of steps. From a model selection point of view, the D-EM fails
to eliminate the redundant components while the PD algorithm provides
a parsimonious mixture fit.

\subsection{Simulation Experiment}
\label{sec-sim}

We consider the simulated data described in Section
\ref{sec-criteria}. The data were generated from the multivariate
normal mixture densities $N_p(\mQ, {\bf I})$ with $p = 3$ and $n =
270$ (see Table \ref{table:compare}) by selecting the true mixing
measure for $\mQ \in \mG$ as the coordinates of $\boldsymbol{\theta}_j
\, (j = 1, \ldots, m)$ from the set $\{-5, 0, 5\}$ in all possible
combinations, with equal mass at each support vector.

Following Section \ref{sec-tree}, we apply the Penalized Dual
algorithm in the context of building a collection of semiparametric
mixture models for selected values of the sieve parameter
$\delta$. This will yield estimators with both many and few active
support vectors; thereby providing a mechanism to demonstrate the
superiority of our method over the D-EM algorithm across a range of
applications. Both the PD and PD$^{\mathrm{IC}}$ algorithms provide
uniformly better performance, producing 6 to 40-fold improvement in
CPU factor over the D-EM algorithm. As illustrated in Section
\ref{sec-tree}, in fitting $N_p(\mQ, \, \delta \, {\boldsymbol
\Sigma})$, there is a trade-off between decrease in the sieve
parameter $\delta$ and the increase in mixture complexity $m$; by
increasing $\delta$, we obtain a reduction in the mixture complexity
$m$. 
\begin{table}[htp]
\caption{Building a collection of semiparametric mixture models
$\left\{ \left(\wQ_{\delta}, \, \delta \, \wS \right)\!: \delta \in
\Re_{+} \right\}$ using Algorithm 1 for the simulated, Fisher iris and
Yeast microarray data sets. First step involved choosing the observed
data matrix ${\bf y}$ for $\boldsymbol{\Theta}_m$ and setting $\wS =
{\bf S}$ in finding $l(\wmbp)$ using the PD and discrete EM (D-EM)
algorithms. The estimate $\wmbp$ obtained from the PD algorithm with
the corresponding fixed support set and $\delta \, {\bf S}$ are
employed as parameter starting values for finding
$l\left(\wQ_{\delta}, \delta \, \wS \right)$ using the continuous EM
(C-EM) algorithm.} 
\label{table:compare1} \vspace{.3cm}
\centerline{\small
\begin{tabular}{llccccc} \hline 
{\bf Data} & {\bf Algorithm} & \multicolumn5c{$\delta$} \\
& & 5 & 2 & 1 & 0.5 & 0.2 \\ \hline
Simulated & PD & -2642.8555  &  -2393.6817 & -2313.6291 & -2278.7175 &
-2178.5765 \\
& D-EM & -2642.8604 & -2393.6822 & -2313.6299 & -2278.7175 & -2178.5766
\\
& C-EM & -2313.2 & -2313.2 & -2192.13 & -2053.37 & -1876.04 \\
Fisher Iris & PD & -629.1448 & -449.8594 & -376.9440 & -311.5519 &
-192.0285 \\  
  & D-EM & -629.1496 & -449.8595 & -376.9442 & -311.5520 & -192.0285 \\
  & C-EM & -379.91 & -217.3 & -149.63 & -49.16  & -136.65 \\ 
Yeast Microarray & PD & -8088.9982 & -5371.8998 & -3691.6696 &
  -1798.2265 & - \\
                 & D-EM & -8088.9987 & -5371.8999 & -3691.6696 &
		 -1798.2265 & - \\
		 & C-EM & -4025.3 & -2626.1 & -142.2 & 6544.0 & - \\
		 \hline 
\end{tabular} }
\end{table}

As discussed in Section \ref{sec-criteria}, an important attribute
of the convergence of an algorithm is the value of loglikelihood, as it
indicates accuracy on a confidence interval scale. Therefore, in Table
\ref{table:compare1}, we present the loglikelihood values obtained 
using various algorithms. The CPU factor for the PD algorithm over the
EM algorithm ranged from ten to over forty-fold for $\delta
\in \{5, 2, 1, 0.5, 0.2\}$. As predicted by the theory, for the PD and
D-EM algorithms, the final $\Psi(\mQ^{(t)})$ given by (\ref{eq:gstop})
does provide a guarantee on the level of algorithmic convergence
$\Lambda^{(t)}$. Indeed, in some cases the bound $\Lambda^{(t)} \leq
\Psi(\mQ^{(t)})$ was very conservative. Moreover, the convergence
criteria for the PD algorithm described in Section \ref{sec-pdalgor}
achieved the desired accuracy in $\Lambda^{(t)}$; however, typically
the PD algorithms terminated at a considerably higher accuracy than
the D-EM algorithm.  In order to measure this effect, we continued the
D-EM algorithm to the same level of accuracy as that of the PD for
$\delta = 1$. In this case, for the EM algorithm, N$^{(t)} =
9,987$ at convergence, resulting in a CPU factor of 60 instead of 24. 

\subsection{Fisher Iris and Yeast Microarray Data Sets}
\label{sec-iris}

We fit a mixture of multivariate normal distributions to the Fisher
iris data described earlier by finding the NPMLE of $\mQ$ and by
building a collection of semiparametric mixture models for selected
values of $\delta$; results are presented in Table
\ref{table:compare1}. The performance of the PD and D-EM algorithms
was similar to that of the simulated data.

Instead of choosing the PD solution as the parameter starting values
for the C-EM, we consider random values to demonstrate their effect on
a given algorithm. It is important to recognize that the C-EM
algorithm requires an a priori knowledge of $m$. We considered the
Fisher iris data with $\delta = 1$ and randomly selected $m = 15$ data
vectors from $n = 150$ as parameter starting values for the algorithm.
In the ten runs of the C-EM algorithm with random starting values,
$l(\wQ, \wS)$ ranged from -151.99 to -179.13; all of which are
sub-optimal modes due to the ``poor choice'' of starting
values. Similar behavior of the C-EM algorithm, in reaching a
sub-optimal solution, was observed by \cite{pilla:01} for the galaxy
data. Without an a priori knowledge of $m$, choosing $m$ can be quite
a challenge for using the C-EM algorithm in large-scale practical
problems.

The main technology for conducting high-throughput experiments in
functional genomics is the microarray---a technical approach for
assaying the abundance of mRNA for several genes simultaneously (see
\cite{hastie:01} for literature). A gene expression data set collects
the expression values from a series of DNA microarray experiments with
each column representing an experiment. Analysis of the expression
patterns obtained from large gene arrays reveal the existence of
clusters of genes with similar expression patterns. It is common to
write the gene expression data of $n$ genes, each measured at $k$
individual array experiments (e.g., single time points or conditions)
as an $n \times k$ matrix. \cite{holter:00} analyzed a subset of the
original published yeast cdc15 cell-cycle data which consist of $n =
696$ genes under $p = 12$ time points or conditions. An important
scientific question is to find out which genes are most similar to
each other, in terms of their expression profiles across samples. One
way to organize gene expression data is to cluster genes on the basis
of their expression patterns.  One can think of the genes as points in
$\Re^{12}$, which we want to cluster together in some fashion.

We fit multivariate normal mixtures to the yeast microarray data by
finding the NPMLE of $\mQ$. This is an example of high-dimensional
modeling. We observed similar performance of the algorithms to the
previous examples. Table \ref{table:compare1} presents the
loglikelihood values. Since each observation is a point in $\Re^{12}$,
at $\delta = 0.2$, we obtain the empirical CDF as the MLE. That is,
each observation is its own component; hence the solution is not
interesting.

\section{Discussion}
\label{sec-discuss}

In this article we developed a framework for approximating the
continuous parameter space and created an algorithm (based on the
Penalized Dual method) for finding the maximum of $l(\mQ)$;
consequently an algorithm for estimating the mixture complexity. We
established convergence properties of the proposed algorithm. By
exploiting the inherent advantage of the penalty formulation, we
derived a technique for converting the parameter estimators from the
Penalized Dual problem into those for the mixture probability
parameters. We established the existence of parameter estimators and
derived convergence results for the Penalized Dual algorithm, for
fitting overparameterized mixture models.  It was shown empirically
that the Penalized Dual algorithm has a faster rate of convergence,
compared with the discrete EM algorithm for overparameterized mixture
problems.

The algorithm based on the Penalized Dual method reaches closer to the
global maximum and is robust to the choice of the support set
${\boldsymbol \Theta}_m \subset \Omega$ (dimensionality of the
problem). These are desirable features for (1) analyzing
high-dimensional data, and (2) for building a collection of
semiparametric mixture models. The dimension of the dual optimization
problem is fixed at $d$, the number of distinct observed data vectors;
whereas that of the discrete EM grows with the cardinality of
${\boldsymbol \Theta}_{m}$. For discrete mixture problems, such as
binomial or Poisson, often $d \ll n$; therefore, there is no
dimensionality cost with the dual problem. When the cardinality of
${\boldsymbol \Theta}_{m}$ is large, the discrete EM algorithm fails
to converge to the MLE, for all practical purposes, in certain mixture
problems. 

We derived several important structural properties of multivariate
normal mixtures in which $\mQ$ is modeled nonparametrically in the
presence of an unknown variance-covariance matrix ${\boldsymbol
\Sigma} \in \mS$ common to all $m$ components. The role of the sieve
parameter in reducing the dimension of the mixture problem was
demonstrated by creating new graphical devices, namely the Mixture
Tree plots.

The proposed methods are very powerful in searching over the whole
discretized parameter space and in yielding a parsimonious mixture
model. The discrete EM algorithm can be very difficult, if not
impossible, in yielding a parsimonious model in problems with hundreds
or thousands of parameters. Such problems are becoming increasingly
common due to the rapid explosion of high-throughput data in
microarray data and data mining. The applications for the methods
described in this article are rich. Multivariate normal mixtures arise
in many different practical scenarios, including data mining,
knowledge discovery, data compression, pattern recognition and pattern
classification. 

\newpage
\centerline{\large \bf Appendix: Technical Derivations}
\setcounter{equation}{0}\def\theequation{A.\arabic{equation}}
\vspace{0.1in}

\noindent{\bf A.1 Relation Between the Primal and Dual Problems}

We establish the relation between the primal and dual problems at the
solution using the change of variable ${\rm g}_{_\mQ}({\bf y}_i) =
(n_i/n) (w_i)^{-1}$ $(i = 1, \ldots, d)$. As a first step, we prove the
following claim.

{\em Claim.} The maximization of the primal problem in
(\ref{eq:primal}) is equivalent to
\begin{eqnarray}  
   \label{eq:dual2}
   \underset{{\bf g}_{_\mQ}}{\mbox{min}} \, \sum_{i = 1}^d n_i \;
   \ln \, \{{\rm g}_{_\mQ}({\bf y}_i) \}
\end{eqnarray}
subject to ${\bf g}_{_\mQ} = \{{\rm g}_{_\mQ}({\bf y}_1), \ldots,
{\rm g}_{_\mQ}({\bf y}_d)\}^T \in \Re_+^d$ and
$\mD_{\mQ}(\boldsymbol{\theta}_j) \leq 0$ for $\boldsymbol{\theta}_j
\in {\boldsymbol \Theta}_m$ ($j = 1, \ldots, m$), where 
$\mD_{\mQ}(\boldsymbol{\theta}_j)$ is defined in (\ref{eq:grad}).

The gradient constraints $\mD_{\mQ}(\boldsymbol{\theta}_j) \leq 0$ can
equivalently be expressed as
\begin{eqnarray}
   \label{eq:constr}
   \sum_{i = 1}^d \left(\frac{n_i}{n} \right)
   \frac{f_{\boldsymbol{\theta}_j}({\bf y}_i)}{{\rm g}_{_\mQ}({\bf
   y}_i)} \leq 1 \quad \mbox{for}  \quad j = 1, \ldots, m. 
\end{eqnarray}

Let ${\mQ}^{\star} \in \mG$ be the solution to the primal problem in
(\ref{eq:primal}) and let $\mQ \in \mG$ be any solution that satisfies
constraints of the dual problem in (\ref{eq:dual2}). The equivalence
between the primal problem in (\ref{eq:primal}) and the dual problem
in (\ref{eq:dual2}) follows by establishing that
\begin{eqnarray}
   \label{eq:ineql}
   \sum_{i = 1}^d n_i \; \ln \, \{{\rm g}_{_\mQ}({\bf  y}_i) \} \geq
   \sum_{i = 1}^d n_i \; \ln \, \{{\rm g}_{_{\mQ^{\star}}}({\bf  y}_i)
   \}. 
\end{eqnarray}
Since $\ln \, (x + \lambda) \geq \ln \, x + \lambda/(x + \lambda)$,
where $x + \lambda = {\rm g}_{_\mQ}$ and $x = {\rm
g}_{_{\mQ^{\star}}}$, the above inequality yields 
\begin{eqnarray}  
\label{eq:ineqal} 
  \sum_{i = 1}^d n_i \; \ln \, \{{\rm g}_{_\mQ}({\bf y}_i)\} &\geq&
  \sum_{i = 1}^d n_i \; \ln \, \{{\rm g}_{_{\mQ^{\star}}}({\bf y}_i)
  \} +  \sum_{i = 1}^d n_i \; \frac{ \{{\rm g}_{_\mQ}({\bf y}_i) -
  {\rm g}_{_{\mQ^{\star}}}({\bf y}_i) \}}{{\rm g}_{_\mQ}({\bf y}_i)}
  \notag \\  &=& \sum_i n_i \; \ln \, \{{\rm g}_{_{\mQ^{\star}}}({\bf
  y}_i) \} - \sum_i n_i \; \left\{ \frac{\sum_j \pi_j \,
  f_{\boldsymbol{\theta}_j}({\bf y}_i)}{{\rm g}_{_\mQ}({\bf y}_i)} - 1
  \right\}.  
\end{eqnarray}
The second term in the right-hand side of (\ref{eq:ineqal})
is less than zero since $\mD_{\mQ}(\boldsymbol{\theta}_j) \leq 0$ and
hence the relation (\ref{eq:ineql}) holds. Therefore, the claim is
established. 

Define $w_i = (n_i/n) \{{\rm g}_{_\mQ}({\bf y}_i)\}^{-1}$ so that the
constraints in (\ref{eq:constr}) become $\sum_i w_i \,
f_{\boldsymbol{\theta}_j}({\bf y}_i) \leq 1$ for $j = 1, \ldots,
m$. From this definition of $w_i$, the dual problem in
(\ref{eq:dual2}) can be expressed as 
\begin{eqnarray*}
   \underset{{\bf w} \, \in \, \Re^d_+}{\mbox{min}} \, \left\{ \sum_{i
   = 1}^d n_i \, \ln \, \left( \frac{n_i}{n} \right) - \sum_{i = 1}^d
   n_i \, \ln \, (w_i) \right\}.
\end{eqnarray*}
Equivalently, the problem is $\underset{{\bf w}}{\mbox{max}} \,
\sum_i n_i \, \ln \, (w_i)$ subject to ${\bf w} \in \Re^d_+$ which is
the dual optimization problem in (\ref{eq:dual}).

\vspace{0.15in}
\noindent{\bf A.2 Derivation of the Penalized-Dual Estimator
$\widehat{\pi}^{\star}_{j, \gamma}$}    

First, from the primal-dual relationship, it follows that 
\begin{eqnarray}
  \label{eq:pdr}
  \widehat{w}_i = \frac{n_i}{n} \, \left[ \sum_{j = 1}^m
  \widehat{\pi}_j \, \left\{f_{{\boldsymbol \theta}_j}({\bf y}_i)
  \right\} \right]^{-1} \quad \mbox{for} \quad {\bf y}_i \in \mY; \, i
  = 1, \ldots, d. 
\end{eqnarray}

Second, the following fixed-point equation is obtained by solving
(\ref{eq:deriv1}),  
\begin{eqnarray}  
  \label{eq:fixp}
  \widehat{w}_{i, \gamma} = \frac{n_i}{n} \left[ \sum_{j = 1}^m
  \left\{p_{_{\boldsymbol{\theta}_{j}}} \left(\widehat{\bf
  w}_{\gamma} \right) \right\}^{\, (\gamma - 1)} \, f_{{\boldsymbol
  \theta}_j}({\bf y}_i) \right]^{-1} \quad \mbox{for} \quad i = 1,
  \ldots, d.  
\end{eqnarray}
By comparing the right-hand sides of (\ref{eq:pdr}) and
(\ref{eq:fixp}), it is clear that ${\rm g}_{_{\wQ}}({\bf y}_i)$
parallels the term $\sum_{j} \{
p_{_{\boldsymbol{\theta}_{j}}}(\widehat{\bf w}_{\gamma})\}^{\, (\gamma
- 1)} \; f_{\boldsymbol{\theta}_j}({\bf y}_i)$ and that the latter
expression resembles a mixture density with
$\{p_{_{\boldsymbol{\theta}_{j}}}(\widehat{\bf w}_{\gamma}) \}^{\,
(\gamma - 1)}$ playing the role of $\widehat{\pi}_j$. 

\vspace{0.15in}
\noindent{\bf A.3 Hessian Matrix of the Function $\mK({\bf z}, \gamma)$}

Let ${\bf F} = ({\bf f}_{\boldsymbol{\theta}_{1}}, \ldots, {\bf
f}_{\boldsymbol{\theta}_{m}})^{T}$ be an $(m \times d)$ matrix where
${\bf f}_{\boldsymbol{\theta}_{j}} =
\{f_{\boldsymbol{\theta}_{j}}({\bf y}_{1}), \ldots,
f_{\boldsymbol{\theta}_{j}}({\bf y}_{d})\}^{T}$ is the $d$-dimensional
vector. In the sequel, we denote a vector of ones by ${\bf 1}$ (with
dimension clear from the context) and the diagonal matrix with
elements ${\bf x}$ by $\diag({\bf x})$. Therefore, 
\begin{equation*} 
   \mK({\bf z}, \gamma) = \frac{1}{n} {\bf n}^{T} \cdot {\bf z} -
   \frac{1}{\gamma} \; {\bf 1}^{T} \cdot {\bf p}^{\, \gamma} \quad
   \mbox{for} \quad {\bf z} \in \Re \; \mbox{and} \; \gamma \in
   \Re_+, 
\end{equation*}
where ${\bf n} = (n_{1}, \ldots, n_{d})^{T}$ and the constraint vector
\begin{equation*}
  {\bf p} = \left\{ p_{_{\boldsymbol{\theta}_{1}}}({\bf z}), \ldots,
  p_{_{\boldsymbol{\theta}_{m}}}({\bf z})\right\}^{T} 
\end{equation*}
with $p_{_{\boldsymbol{\theta}_{j}}}({\bf z})$ (sometimes
written as $p_{j}$ for exposition) is as in (\ref{eq:pgrad}) expressed
in terms of ${\bf z} \in \Re$. 

The Hessian matrix of $\mK({\bf z}, \gamma)$ has the following
elements:  
\begin{eqnarray}  
   \frac{\partial}{\partial {\bf z}} \mK({\bf z}, \gamma) &=&
   \frac{{\bf n}}{n} - \diag({\bf w}) \cdot 
   \left\{ {\bf F}^{T} \cdot {\bf p}^{\, (\gamma
   -1)} \right\}, \notag \\ 
   \frac{\partial^2}{\partial {\bf z} \; \partial {\bf z}^{T}}
   \mK({\bf z}, \gamma) &=& - \; \diag({\bf w} \cdot \left\{{\bf
   F}^{T} \cdot {\bf p}^{\, (\gamma -1)} \right\} \notag \\ 
   && -(\gamma - 1) \; \diag({\bf w}) \cdot
   {\bf F}^{T} \; \diag \left\{ {\bf p}^{\,
   (\gamma - 2)} \right\} \; {\bf F} \cdot
   \diag({\bf w}), \notag \\
   \frac{\partial}{\partial \gamma} \mK({\bf z}, \gamma) &=&
   \frac{1}{\gamma} \sum_{j = 1}^m (p_j)^{\, \gamma}
   \left\{ \frac{1}{\gamma} - \ln \, (p_j) \right\},  \notag \\
   \frac{\partial^2}{\partial \gamma^2} \mK({\bf z}, \gamma) &=&
   \frac{1}{\gamma} \left[ - \frac{\partial \mK({\bf z},
   \gamma)}{\partial \gamma} + \sum_{j = 1}^m (p_j)^{\, \gamma}  
   \left\{ \frac{1}{\gamma} - \ln \, (p_j) \right\} \ln \, (p_j) -
   \frac{1}{\gamma^2} \right], \notag 
\end{eqnarray}
where ${\bf w} \in \Re^d_+$ is expressed as $\{\exp({\bf z}_1),
\ldots, \exp({\bf z}_d)\}$ and 
\begin{eqnarray}  
  \label{eq:hessian}
   \frac{\partial^2}{\partial {\bf z} \; \partial \gamma} \mK({\bf z},
   \gamma) &=& - \frac{1}{\gamma} \sum_{j = 1}^m (p_j)^{\, (\gamma -
   1)} \; \ln \, (p_j). 
\end{eqnarray}

\vspace{0.15in}
\noindent{\bf A.4 Proofs}

In the sequel, we write $p_{_{\boldsymbol{\theta}_{j}}}(\widehat{\bf
w}_{\gamma}) = \widehat{p}_{j, \gamma}$ for 
exposition. 

{\em Proof of Theorem \ref{lemm1}.} From the fixed-point equation
(\ref{eq:fixp}), we have the EM solution 
\begin{eqnarray}
   \widehat{\pi}_{_{j, \mathrm{EM}}} &=& \widehat{\pi}_{j,
   \gamma}^{\dag} \cdot \sum_{i = 1}^d \left(\frac{n_i}{n} \right)
   \frac{f_{\boldsymbol{\theta}_j}({\bf y}_i)}{\sum_{k = 1}^m 
   \, \widehat{\pi}_{k, \gamma}^{\dag} \;
   f_{\boldsymbol{\theta}_k}({\bf y}_i)} \quad \mbox{for} \quad
   \boldsymbol{\theta}_j \in {\boldsymbol \Theta}_m \, (j = 1,
   \ldots, m) \nonumber \\ 
   &=& \wp_{\gamma} \; \{\widehat{p}_{j, \gamma}\}^{\, (\gamma - 1)}
   \sum_i \left(\frac{n_i}{n} \right) \frac{\wp_{\gamma}^{-1} \,
   f_{\boldsymbol{\theta}_j}({\bf y}_i)}{\sum_k 
   \{ \widehat{p}_{k, \gamma} \}^{\, (\gamma - 1)} \,
   f_{\boldsymbol{\theta}_k}({\bf y}_i)} 
\end{eqnarray}
which follows from (\ref{eq:ems}), where $\wp_{\gamma}$ is given in
(\ref{eq:sum}). From (\ref{eq:fixp}), the last equation becomes 
\begin{eqnarray*}
  \{\widehat{p}_{j, \gamma}\}^{\, (\gamma - 1)} \sum_{i = 1}^d w_i \, 
  f_{\boldsymbol{\theta}_j}({\bf y}_i).
\end{eqnarray*}
This again simplifies to $\{\widehat{p}_{j, \gamma} \}^{\,
(\gamma - 1)} \; \widehat{p}_{j, \gamma} = \{\widehat{p}_{j, \gamma}
\}^{\,\gamma}$ due to the relationship in (\ref{eq:pgrade}). Thus
$\widehat{\pi}_{_{j, \mathrm{EM}}} =  \{ \widehat{p}_{j, \gamma}\}^{\,
\gamma} = \widehat{\pi}_{j, \gamma}^{\star}$ and the proof of part (a)
follows. As a consequence of the EM result, the estimators are in the
unit simplex ${\boldsymbol \Pi}^{\star}$ as claimed in part (b). Proof
of part (c) follows by using the first inequality in part (b) in
conjunction with relation (\ref{eq:pdest}). These two imply that
$p_{_{\boldsymbol{\theta}}} (\widehat{w}_{i, \gamma}) \leq 1$ for all
$\boldsymbol{\theta} \in {\boldsymbol \Theta}_m$. Hence, our estimator
is in the feasible region as claimed. \hfill
\rule{3mm}{3mm}  

\vspace{0.1in}
Next, we need the following lemma to prove Theorem
\ref{thm-result}. 
\begin{lemma}
\label{lem:lresult} 
{\rm As $\gamma \rightarrow \infty$, $\wp_{\gamma}^{-1}\rightarrow 1$.
}
\end{lemma}
{\em Proof.} From the following Lyapunov's inequality
\citep{lehmann:99},  
\begin{eqnarray*}
   E \Big(X^{(\gamma - 1)} \Big)^{\frac{1}{(\gamma - 1)}} \leq E
   \Big(X^{\gamma} \Big)^{\frac{1}{\gamma}},
\end{eqnarray*}
one can find a bound for $\wp_{\gamma}$. For $m$ number of
constraints, it follows that 
\begin{eqnarray*}
   \left(\sum_{j = 1}^m \frac{1}{m} \; \{\widehat{p}_{j, \gamma}
   \}^{\, (\gamma - 1)} \right)^{\frac{1}{(\gamma - 1)}} 
   &\leq& \left( \sum_{j = 1}^m \frac{1}{m} \; \left\{ \widehat{p}_{j,
   \gamma} \right\}^{\, \gamma} \right)^{\frac{1}{\gamma}} =
   m^{-\frac{1}{\gamma}} \cdot 1. 
\end{eqnarray*}
Equivalently, $\sum_{j} \{\widehat{p}_{j, \gamma} \}^{\,
(\gamma - 1)} \leq m^{\frac{1}{\gamma}}$. Hence as $\gamma \rightarrow
\infty$, we obtain $\wp_{\gamma}^{-1} \rightarrow 1$. \hfill
\rule{3mm}{3mm} 

\vspace{0.15in} 
{\em Proof of Theorem \ref{thm-result}.} From Corollary
\ref{cor:mgrad}, we have   
\begin{eqnarray*}
   \mD_{\widehat{\mQ}^{\dag}_{\gamma}}(\boldsymbol{\theta}_j) \leq
   \wp_{\gamma} - 1 \quad \mbox{for} \quad \boldsymbol{\theta}_j \in
  {\boldsymbol \Theta}_m \, (j = 1, \ldots, m),
\end{eqnarray*}
where $\mD_{\widehat{\mQ}^{\dag}_{\gamma}}(\boldsymbol{\theta}_j)$
and $\wp_{\gamma}$ are defined in (\ref{eq:modgrad}) and
(\ref{eq:sum}), respectively.  

From Lemma \ref{lem:lresult}, we have $\wp_{\gamma}^{-1} \rightarrow
1$ as $\gamma \rightarrow \infty$. Therefore, in the limit, the
primal-gradient function satisfies the inequality  
\begin{eqnarray*}
  \underset{\gamma \rightarrow \infty}{\mathrm{lim}}
  \mD_{\widehat{\mQ}^{\dag}_{\gamma}}(\boldsymbol{\theta}_j) \leq 0
  \quad \mbox{for} \quad \boldsymbol{\theta}_j \in
  {\boldsymbol \Theta}_m \, (j = 1, \ldots, m).
\end{eqnarray*}
The compactness of the parameter space $\boldsymbol{\Pi}$ can in turn
be used to establish the convergence of ${\bf
g}_{_{\wQ^{\dag}_{\gamma}}}$ to the maximizing value ${\bf
g}_{_{\wQ}}$. If the vector of masses $\mbp$ for ${\bf g}_{_{\wQ}}$
are uniquely determined, then the masses must converge as well. This
in turn implies that as $\gamma \rightarrow \infty$, the mixing
distribution $\widehat{\mQ}^{\dag}_{\gamma}$ with
$\widehat{\boldsymbol{\pi}}^{\dag}$ as the vector of masses is the
NPMLE. Consequently, $\widehat{\pi}_{j, \gamma}^{\dag} \rightarrow
\widehat{\pi}_j$ as $\gamma \rightarrow \infty$ for $j = 1, \ldots,
m$. \hfill \rule{3mm}{3mm}    

\vspace{0.15in} 
{\em Proof of Theorem \ref{thm-grad}.} From (\ref{eq:modgrad}) and the
relation $\widehat{\pi}_{j, \gamma}^{\dag} = \{\widehat{p}_{j, \gamma}
\}^{\, (\gamma - 1)} \wp_{\gamma}$, it follows that  
\begin{eqnarray*}
  \mD_{\widehat{\mQ}^{\dag}_{\gamma}}(\boldsymbol{\theta}_j) =
  \sum_{i = 1}^d n_i \left[ \frac{\wp_{\gamma}^{-1} \, 
  f_{\boldsymbol{\theta}_j}({\bf y}_i)}{\sum_k \left\{\widehat{p}_{k,
  \gamma} \right\}^{\, (\gamma - 1)} \, f_{\boldsymbol{\theta}_k}({\bf
  y}_i)} - 1 \right] \quad \mbox{for} \quad \boldsymbol{\theta}_j \in
  {\boldsymbol \Theta}_m. 
\end{eqnarray*}

From (\ref{eq:fixp}), we have $\widehat{w}_{i, \gamma} = n_i
\left[\sum_k \left\{ \widehat{p}_{k, \gamma} \right\}^{\, (\gamma -
1)} f_{\boldsymbol{\theta}_k}({\bf  y}_i) \right]^{-1}$ and hence the
last equation simplifies to  
\begin{eqnarray*}
   \wp_{\gamma}^{-1} \sum_{i = 1}^d \widehat{w}_{i, \gamma} \,
   f_{\boldsymbol{\theta}_j}({\bf y}_i) - 1 \quad \mbox{for} \quad
   \boldsymbol{\theta}_j \in {\boldsymbol \Theta}_m. 
\end{eqnarray*}
The desired result follows from the definition of $\widehat{p}_{j,
\gamma}$. \hfill \rule{3mm}{3mm}  

\vspace{0.15in} 
{\em Proof of Theorem \ref{lem-concave}.} 
The proof of part (a) is a consequence of the negative definiteness of
${\bf H}$ and the properties of $\mP({\bf w}, \gamma)$ given in
Proposition \ref{lem-penalty}. For part (b), let 
\begin{eqnarray*}
   \widehat{{\bf z}}_{\gamma} = \underset{{\bf
   z} \in \Re}{\mathrm{arg \; max}} \; \mK({\bf z}, \gamma)
\end{eqnarray*}
for any fixed $\gamma \in \Re_+$. That is, $\widehat{{\bf
z}}_{\gamma}$ is the maximizer of the $\mK({\bf z}, \gamma)$ for a
fixed $\gamma$. From equation (\ref{eq:kfunc}) for a given
$\widehat{{\bf z}}_{\gamma}$, it follows that
\begin{eqnarray*}
   \frac{\partial}{\partial \gamma} \mK \left(\widehat{{\bf
   z}}_{\gamma}, \gamma \right) = \frac{1}{\gamma} \sum_{j = 1}^m
   \left\{\widehat{p}_{j, \gamma} \right\}^{\gamma} \left\{
   \frac{1}{\gamma} - \ln \, \left(\widehat{p}_{j, \gamma} \right)
   \right\},  
\end{eqnarray*}
where $\widehat{p}_{j, \gamma}$ is expressed in terms of
$\widehat{{\bf z}}_{\gamma}$. From part (c) of Theorem \ref{lemm1}, we
have $\{\widehat{p}_{j, \gamma}\}^{\, \gamma} \leq 1$ for $j = 
1, \ldots, m$. Therefore, $\ln \, (\widehat{p}_{j, \gamma}) \leq
0$ and hence $\partial \mK(\widehat{{\bf z}}_{\gamma},
\gamma)/\partial \gamma > 0$ for any finite $\gamma \in \Re_+$ and
fixed $\widehat{{\bf z}}_{\gamma}$. That is, for a fixed ${\bf z} \in
\Re$, the only point at which the function $\mK({\bf z},
\gamma)$ can approach its supremum is at $\gamma = \infty$. \hfill  
\rule{3mm}{3mm}  

\vspace{0.15in} 
{\em Proof of Theorem \ref{thm-mmml}.} For exposition we drop the
subscript $\bS$ from $\mQ^{\star}$ and $\wQ$; however it is understood
that the mixing measures have a dependence on the fixed $\bS$.

Part (1): We first establish that $\wQ$ is an NPMLE. We start with
creating a path in $\mG$ from $\wQ$ to $\mQ^{\star}$, by letting
$\mQ_{\alpha} = (1 - \alpha) \, \wQ + \alpha \, \mQ^{\star}$ for
$\alpha \in [0, 1], \mQ^{\star} \in \mG$ and $\wQ$ satisfying the
relation (\ref{gradineq}). Note that $\mQ_{\alpha} \in \mG$;
therefore, $\mG$ is a convex set. The loglikelihood along this path
satisfies
\bes
  l \left(\mQ_{\alpha}; \bS \right) \geq (1 - \alpha) \, l \left(\wQ;
  \bS \right) + \alpha \, l \left(\mQ^{\star};\bS \right) 
\ees
for $\alpha \in [0, 1]$ and for a fixed $\bS$. Therefore,
$l(\mQ_{\alpha}; \bS)$ is a concave function for a fixed $\bS$. The
{\em directional derivative} of $l(\mQ_{\alpha}; \bS)$ at ${\rm g}_{_{\wQ}}$
toward ${\rm g}_{_{\mQ^{\star}}}$ can be expressed, after
simplification, as 
\bea
  \frac{d}{d \alpha} \, l(\mQ_{\alpha}; \bS)\Big|_{\alpha = 0} &=&
  \label{eq:int}
  \int \mD_{\widehat{\mQ}}(\boldsymbol{\mu}; {\boldsymbol \Sigma})
  \; d \, \mQ^{\star}(\boldsymbol{\mu}).
\eea
From (\ref{gradineq}), it follows that (\ref{eq:int}) is $\leq 0$ for
all $\mQ^{\star} \in \mG$ and $\wQ$ satisfying (\ref{gradineq}). This
result combined with the concavity of $l(\mQ_{\alpha}; \bS)$ implies
that $\wQ$ is an NPMLE of $\mQ \in \mG$.

Part (2): We establish the uniqueness of the NPMLE. Suppose $\wQ$ and
$\mQ^{\star}$ are two NPMLEs of $\mQ \in \mG$, then
\bes
 l \left(\mQ_{\alpha}; \bS \right) = (1 - \alpha) \, l \left(\wQ;
 \bS \right) + \alpha \, l \left(\mQ^{\star};\bS \right) 
\ees
for all $\alpha \in [0, 1]$ and for a fixed $\bS$. This implies that 
the derivative $d \, l(\mQ_{\alpha}; \bS)/d \, \alpha$ at
$\alpha = 0$ is exactly zero. This implies that $\mQ^{\star}$ (and
$\wQ$) is supported on $\{\boldsymbol{\mu}_1, \ldots,
\boldsymbol{\mu}_K\}$ [i.e., the set of zeroes of
$\mD_{\widehat{\mQ}}(\boldsymbol{\mu}; {\boldsymbol
\Sigma})$]. Furthermore, the second derivative  
\bes
  \label{eq:dsq}
  \frac{d^2}{d \, \alpha^2} \, l(\mQ_{\alpha}; \bS)\Big|_{\alpha = 0}
  &=& 0,
\ees
which implies that 
\bes
  - \sum_{i = 1}^n \frac{ \left\{ {\rm g}_{_{\mQ^{\star}}}({\bf
  y}_i; \bS) - {\rm g}_{_{\wQ}}({\bf y}_i; \bS) \right\}^2}{ \left\{ 
  \alpha \, {\rm g}_{_{\mQ^{\star}}}({\bf y}_i; \bS) + (1 - \alpha) \,  
  {\rm g}_{_{\wQ}}({\bf y}_i; \bS) \right\}^2 } = 0. 
\ees
That is, 
\bes
  \label{eq:gdenin}
  {\rm g}_{_{\mQ^{\star}}}({\bf y}_i; \bS) = {\rm g}_{_{\wQ}}({\bf
  y}_i; \bS) \quad \mbox{for all} \quad i = 1, \ldots, n.
\ees
However, the linear independence of the vectors ${\bf
f}_{\boldsymbol{\mu}_j}({\bf y}; {\boldsymbol \Sigma})$ for $j = 1,
\ldots, K$ implies that \\ $\{\pi_1, \ldots, \pi_K\} =
\{\pi_1^{\star}, \ldots, \pi_K^{\star}\}$. That is, $\mQ^{\star} =
\wQ$; establishing the uniqueness of the NPMLE of $\mQ$ for a
fixed $\bS$.
\hfill \rule{3mm}{3mm}  

\vspace{0.2in}
\bibliographystyle{apalike}
\bibliography{pillaetal}

\end{document}